\documentclass[a4paper]{article}

\usepackage{geometry}                
\usepackage{graphicx}
\usepackage{tikz}
\usepackage{enumerate}
\usepackage{amsmath,amsfonts, amssymb,amsthm}
\usepackage{epstopdf}
\usepackage[utf8]{inputenc} 
\usepackage[T1]{fontenc} 
\usepackage[extra]{tipa}
 
\usepackage[english]{babel}

\title{The hyperbolicity of the sphere complex via surgery paths}
\author{Arnaud Hilion\thanks{
Supported by the grant ANR-10-JCJC 01010 of the Agence Nationale de la Recherche.}
, Camille Horbez}

\usepackage[all]{xy}

\usepackage{bbm}
\begin{document}
\maketitle
\newtheorem{de}{Definition} [section]
\newtheorem{theo}[de]{Theorem} 
\newtheorem{prop}[de]{Proposition}
\newtheorem{lemma}[de]{Lemma}
\newtheorem{cor}[de]{Corollary}
\newtheorem{propd}[de]{Proposition-Definition}

\begin{abstract}
In \cite{HM12}, Handel and Mosher have proved that the free splitting complex $\mathcal{FS}_n$ for the free group $F_n$ is Gromov hyperbolic. This is a deep and much sought-after result, since it establishes $\mathcal{FS}_n$ as a good analogue of the curve complex for surfaces.

We give a shorter alternative proof of this theorem, using surgery paths in Hatcher's sphere complex (another model for the free splitting complex), instead of Handel and Mosher's fold paths. As a byproduct, we get that surgery paths are unparameterized quasi-geodesics in the sphere complex.

We explain how to deduce from our proof the hyperbolicity of the free factor complex and the arc complex of a surface with boundary.
\end{abstract}

\section{Introduction}

Let $n\in\mathbb{N}$, and let $M_n=\#_nS^1\times S^2$ be the connected sum of $n$ copies of $S^1\times S^2$.
The fundamental group of the manifold $M_n$ is the free group $F_n$ of rank $n$. A \emph{sphere system} is a collection of disjoint embedded $2$-spheres in $M_n$, none of which bounds a ball, and no two of which are homotopic. In fact, it follows from the work of Laudenbach \cite{Lau73} that two spheres in $M_n=\#_nS^1\times S^2$ are homotopic if and only if they are isotopic:
such spheres are usually called {\em parallel}. More generally, two sphere systems are homotopic if and only if they are isotopic. The \emph{sphere complex} $\mathcal{S}'_n$, introduced by Hatcher in \cite{Hat95}, is the simplicial complex whose $k$-simplices are the homotopy classes of systems of $k+1$ spheres in $M_n$ (a $(k-1)$-dimensional face of a $k$-simplex $\Delta$ is obtained by deleting one sphere of the sphere system corresponding to $\Delta$).

The present paper is devoted to giving a new proof of the following theorem.

\begin{theo}[Handel-Mosher \cite{HM12}]\label{thm:hyperbolicity}
The sphere complex $\mathcal{S}'_n$, equipped with the simplicial metric, is Gromov hyperbolic.
\end{theo}

More precisely, Handel and Mosher \cite{HM12} proved that the free splitting complex $\mathcal{FS}_n$ is Gromov hyperbolic. The free splitting complex $\mathcal{FS}_n$ is the $\text{Out}(F_n)$-complex whose vertices are the nontrivial free splittings of $F_n$ given by a graph of groups decomposition of $F_n$ with a single edge, and two distinct vertices of $\mathcal{FS}_n$ are joined by an edge if the corresponding splittings admit a common refinement (the higher-dimensional simplices are defined in a similar way). However, the complexes $\mathcal{S}'_n$ and $\mathcal{FS}_n$ are 
 canonically
$\text{Out}(F_n)$-equivariantly isomorphic -- see for instance \cite[Lemma 2]{AS11}.  We will denote this isomorphism by $\Theta:\mathcal{S}'_n\rightarrow\mathcal{FS}_n$.
Basically, a nontrivial sphere in $M_n$ gives rise to a one-edge free splitting of $\pi_1(M_n)=F_n$ by van Kampen's theorem, and two nontrivial disjoint nonparallel spheres give rise to two compatible free splittings of $F_n$.

\bigskip

For a long while, there was no fully satisfying metric theory for the group $\text{Out}(F_n)$ of outer automorphisms of $F_n$. It is known that for $n\ge 3$, the group $\text{Out}(F_n)$ is not hyperbolic, and not even relatively hyperbolic -- see \cite{AAS07,BDM09}.
This is one of the main reasons that have led people to look for hyperbolic spaces equipped with a ``good'' action of $\text{Out}(F_n)$, analogous to the complex of curves for mapping class groups \cite{MM99}. The sphere complex $\mathcal{S}'_n$ is a good candidate of a complex with a natural $\text{Out}(F_n)$-action : Aramayona and Souto \cite{AS11} proved that $\text{Out}(F_n)$ is precisely the group of simplicial automorphisms of $\mathcal{S}'_n$.

The geometric study of $\text{Out}(F_n)$ started with the construction by Culler and Vogtmann of an $\text{Out}(F_n)$-complex ``with missing faces'', known as Outer space $CV_n$ \cite{CV86}. There is a deep link between $CV_n$ and $\mathcal{S}'_n$: indeed, the sphere complex $\mathcal{S}'_n$ is the simplicial completion of $CV_n$. 
However, Outer space, endowed with the simplicial metric, is not Gromov hyperbolic: it is possible to construct arbitrarily large flats in $CV_n$. 

In \cite{BF10}, Bestvina and Feighn made an attempt to exhibit hyperbolic $\text{Out}(F_n)$-complexes, however the complexes were in a sense non-canonical. The first striking result was the proof by Bestvina and Feighn of the hyperbolicity of the free factor complex $\mathcal{FF}_n$ \cite{BF11}. Then Handel and Mosher proved Theorem~\ref{thm:hyperbolicity}.
Analysing the natural map from $\mathcal{FS}_n$ to $\mathcal{FF}_n$, Kapovich and Rafi \cite{KR12} have recently shown that the hyperbolicity of $\mathcal{FF}_n$ can be deduced from the hyperbolicity of $\mathcal{FS}_n$.

\bigskip

Both our proof of Theorem~\ref{thm:hyperbolicity} and Handel and Mosher's rely on the criterion of Masur and Minsky \cite{MM99} to prove that a connected simplicial complex $\mathcal{X}$, equipped with the simplicial metric, is Gromov hyperbolic. In order to apply Masur and Minsky's theorem, one has to exhibit a family of ``paths'' in $\mathcal{X}$, and then prove that this family satisfies certain axioms -- see Section~\ref{sec:MM} for details. As an output, one gets that these paths are in fact unparameterized quasi-geodesics in $\mathcal{X}$.

The original proof of Handel and Mosher used ``fold paths'', whereas our proof uses ``surgery paths''. If the former were considered as the ``natural paths'' to move within $\mathcal{FS}_n$, the latter would be the ``natural paths'' in $\mathcal{S}'_n$. A definition of surgery paths is given in Section~\ref{sec-paths}. They have already been shown to be a very useful tool: for instance to prove that the sphere complex is contractible \cite{Hat95}, to prove exponential isoperimetric inequalities for the group $\text{Out}(F_n)$ \cite{HV96},
to study the (un-)distorsion of some natural subgroups of $\text{Out}(F_n)$ \cite{HH11},
and more generally to investigate the geometry of $CV_n$ equipped with the Lipschitz metric \cite{Hor12}.
Fold paths are a slightly delicate elaboration on Stallings' folding paths \cite{Sta83}: their precise definition can be found in Section~2 of \cite{HM12}, and is a bit too long to be recalled here.

As a corollary of Masur and Minsky's theorem, we obtain the following new result:

\begin{theo}\label{thm:uniform}
Surgery paths are unparameterized quasi-geodesics in $\mathcal{S}'_n$ with respect to some uniform constants.  \qed
\end{theo}

Together with Handel and Mosher's results, this implies that there exists a constant $C>0$ (which only depends on $n$) such that for any sphere systems $S_1,S_2\in\mathcal{S}'_n$, corresponding to $T_1=\Theta(S_1),T_2=\Theta(S_2)$ in $\mathcal{FS}_n$, for any surgery path $\gamma$ joining $S_1$ and $S_2$ (in either direction), for any fold path $\gamma'$ joining $T_1$ and $T_2$ (again in either direction), the Hausdorff distance between $\Theta(\gamma)$ and $\gamma'$ is at most $C$. Let us be a bit informal in the comments that follow. Both surgery paths and fold paths can be viewed as sequences of elementary surgeries and elementary folds. The point is that an elementary surgery from a sphere system $S_1$ to a sphere system $S_2$ naturally defines an elementary fold from the tree $T_2=\Theta(S_2)$ to the tree $T_1=\Theta(S_1)$. This indicates that one could think of surgery paths as unfolding paths (i.e. folding paths with reverse orientation), and the striking analogies between our proof of the hyperbolicity of $\mathcal{S}'_n$ and Handel and Mosher's proof of the hyperbolicity of $\mathcal{FS}_n$ consolidates the idea that surgery paths should be unfolding paths in some way. Finding an explicit description of the link between surgery paths and fold paths would be an interesting problem.

\bigskip

One might think at first sight that the realization of the free splitting complex (or, for the matter, of Outer space) as a complex of spheres is loaded with extra geometric data that could overshadow the real geometry of $\text{Out}(F_n)$, already encoded by the graphs (or their corresponding trees) which describe the free splittings of free groups.
However, it turns out that these extra geometric data enable us to shorten the arguments from Handel and Mosher's original proof substantially, and avoid some of the technical difficulties of their proof.

\bigskip

There is a strong correspondence between the important steps of the proof of \cite{HM12} and those of our proof. We have tried to make these correspondences more transparent, in particular by using the same terminology as in \cite{HM12} to name the corresponding phenomena.

\bigskip

The last section of this paper was written during and after the {\em ``Conference on Automorphisms of free groups: Algorithms, Geometry and Dynamics''} that took place at the {\sc crm} in Barcelona in November 2012. It has been deeply influenced by the stimulating atmosphere of the conference. In particular, Section~\ref{sec:arc complex} has been motivated by a question of Karen Vogtmann and comments by Juan Souto and Saul Schleimer. Section~\ref{sec:free factor} has certainly benefited from discussions with Kasra Rafi.

The goal of Section~\ref{sec:other complexes} is to show how to deduce the hyperbolicity of other related complexes from our proof of the hyperbolicity of the sphere complex. The precise definition of these complexes in given in Section~\ref{sec:other complexes}. We recover the result \cite{BF11} of Bestvina and Feighn on the hyperbolicity of the free factor complex, using an approach \`a la Kapovich and Rafi \cite{KR12}, as well as Masur and Schleimer's result \cite{MS13} on the hyperbolicity of the arc complex of a surface with boundary.

\begin{theo} Let $n,g,s$ be positive integers.
\begin{enumerate}[(i)]
  \item {\em (Bestvina-Feighn \cite{BF11})} The free factor complex $\mathcal{FF}_{n}$ is Gromov hyperbolic.
  \item {\em (Masur-Schleimer \cite{MS13})} The arc complex $\mathcal{A}_{g,s}$ (of a surface of genus $g$ with $s$ boundary components) is Gromov hyperbolic.
\end{enumerate}
\end{theo}

\bigskip

{\it Acknowledgements}: We would like to thank the organizers of the {\it Park City Mathematics Institute 2012 Summer Session}. The present work has benefited from discussions we had there: we are grateful to Vincent Guirardel, Lee Mosher, Saul Schleimer and Karen Vogtmann for helpful remarks, and in particular to Ursula Hamenst\"{a}dt for sharing with us her own insight regarding surgery paths and the hyperbolicity of the sphere complex.

\section{The sphere complex}

\subsection{Paths in the sphere complex}

We quickly recall that $M_n=\#_nS^1\times S^2$ is the connected sum of $n$ copies of $S^1\times S^2$. A \emph{sphere system} is a collection of disjoint embedded $2$-spheres in $M_n$, none of which bounds a ball, and no two of which are parallel. The \emph{sphere complex} $\mathcal{S}'_n$ is the simplicial complex whose $k$-simplices are the homotopy classes of systems of $k+1$ spheres in $M_n$. In particular, any $(k-1)$-dimensional face of a $k$-simplex $\Delta$ is obtained by deleting one sphere of the sphere system corresponding to $\Delta$. Throughout the paper, we will sometimes abuse notations and denote in the same way a sphere system $S$ and its homotopy class.

In this paper, we will rather work with the $1$-skeleton $\mathcal{S}_n^{(1)}$ of the barycentric subdivision $\mathcal{S}_n$ of $\mathcal{S}'_n$. Its vertices are the homotopy classes of sphere systems in $M_n$, and two vertices $S$ and $S'$ of $\mathcal{S}_n^{(1)}$ are joined by an edge if the corresponding sphere systems can be homotoped in such a way that one of them is contained in the other (with a slight abuse of notation, we will write $S\subsetneq S'$ or $S'\subsetneq S$). We endow $\mathcal{S}_n^{(1)}$ with the simplicial metric: given two sphere systems $S,S'\in\mathcal{S}_n^{(1)}$, we define $d(S,S')$ to be the length of a shortest path from $S$ to $S'$ in $\mathcal{S}_n^{(1)}$.

\indent A \emph{W-path} in $\mathcal{S}_n^{(1)}$ is a finite sequence of the form $S^0\supseteq S^1\subseteq S^2\supseteq S^3\subseteq S^4$. A \emph{zig-zag path} is a concatenation of W-paths. This notion was intoduced by Handel and Mosher in \cite[section 1.3]{HM12}, where they described geodesics in the free splitting complex in terms of elementary collapses and expansions (our definition of a zig-zag path is coarsely the same as theirs). 

\begin{lemma} \label{ZZ}
Given two sphere systems $S$ and $S'$, there exists a zig-zag path $S=S^0,\dots,S^{2^{k+1}}=S'$ joining $S$ to $S'$ in $\mathcal{S}_n^{(1)}$, where $k\in\mathbb{N}^*$ satisfies $2^{k}-3\le d(S,S')\le 2^{k+1}$.
\end{lemma}

\noindent \textit{Proof}: First note that, given two sphere systems $S$ and $S'$, there exists a zig-zag path joining $S$ to $S'$ in $\mathcal{S}_n^{(1)}$ whose length is at most $d(S,S')+4$. Indeed, whenever there is a path of the form $S^1\subseteq\dots\subseteq S^k$ from $S^1$ to $S^k$ in $\mathcal{S}_n^{(1)}$, then there is also an edge in $\mathcal{S}_n^{(1)}$ joining $S^1$ to $S^k$, so that any geodesic between $S$ and $S'$ has the shape of a zig-zag path, from which we may eventually delete up to $4$ edges at its extremities. We conclude by choosing a zig-zag path joining $S$ to $S'$ of length $l$ which satisfies $d(S,S')\le l \le d(S,S')+4$, choose $k$ such that $2^k<l\le 2^{k+1}$, and if needed, complete the zig-zag path by repeating the last sphere.
\hfill $\square$

\subsection{Hatcher's normal form for sphere systems} \label{sec-NF}

\indent An important tool in the study of sphere systems is Hatcher's normal form (see \cite[Section 1]{Hat95}). Given a sphere system $S$ and a maximal sphere system $\Sigma$, the sphere system $S$ is said to be in \emph{normal form} with respect to $\Sigma$ if every sphere $s\in S$ either belongs to $\Sigma$, or intersects $\Sigma$ transversely in a collection of circles that split $s$ into components called \emph{pieces}, in such a way that for each component $P$ of $M_n\smallsetminus\Sigma$ one has:
\begin{enumerate}[(i)]
\item each piece in $P$ meets each component of $\partial P$ in at most one circle, and

\item no piece in $P$ is a disk which is isotopic, relative to its boundary, to a disk in $\partial P$.
\end{enumerate}

\indent Hatcher proved in \cite[Proposition 1.1]{Hat95} that a sphere system $S$ can always be homotoped into normal form with respect to the maximal sphere system $\Sigma$. Following the exposition in \cite[Section 7.1]{HOP12}, we extend this notion to arbitrary sphere systems in place of $\Sigma$. When $S$, $\Sigma$ are arbitrary sphere systems, we say that $S$ and $\Sigma$ are in normal form if, given any sphere $\widetilde{s}$ in the lift of $S$ to the universal cover of $M_n$ and any sphere $\widetilde{\sigma}$ in the lift of $\Sigma$, the following holds:

\begin{enumerate}[(i)]
\item the spheres $\widetilde{s}$ and $\widetilde{\sigma}$ are either equal, or they intersect in at most one circle, and

\item when $\widetilde{s}$ and $\widetilde{\sigma}$ intersect in one circle, none of the two disks in $\widetilde{s}\smallsetminus\widetilde{\sigma}$ is isotopic, relative to its boundary, to a disk in $\widetilde{\sigma}$.
\end{enumerate}

The argument given in \cite[Lemma 7.2]{HOP12} shows that both definitions agree when $\Sigma$ is a maximal sphere system.
\\
\\ 
\indent Given two sphere systems $S$ and $\Sigma$, we define their \emph{intersection number} $i(S,\Sigma)$ as the minimal number of intersection circles between representatives of the homotopy classes of $S$ and $\Sigma$. Using the same procedure as in the proof of \cite[Proposition 1.1]{Hat95} to remove problematic circles, and an extension of \cite[Proposition 1.2]{Hat95} to the case of arbitrary sphere systems (see the proof of \cite[Lemma 7.3]{HOP12}), one shows that two such representatives minimize the number of intersection circles if and only if they are in normal form (this could also be recovered from \cite{GP09}). Besides, if two sphere systems $S$ and $S'$ are in normal form, then any two of their subsystems will be in normal form too.

\subsection{Surgery paths} \label{sec-paths}

Let $S,\Sigma\in\mathcal{S}_n^{(1)}$, and assume that the sphere system $S$ has been homotoped into normal form with respect to $\Sigma$. As in \cite{HV96}, we now describe a way of passing from $S$ to $\Sigma$ in $\mathcal{S}_n^{(1)}$, using a surgery procedure.  

\begin{figure}[h]
\begin{center}
\begin{tikzpicture}[scale=0.6]
  \draw (0,-0.1) circle (1);
  \filldraw[fill=black!15!white, draw=black] (0,0) circle (0.5);
  \draw (0,0) circle (0.75);
  \draw (2.5,0.2) circle (1);
  \draw (2.5,0.1) circle (0.75);
  \filldraw[fill=black!15!white, draw=black] (2.8,0) circle (0.2);
  \filldraw[fill=black!15!white, draw=black] (2.35,0.15) circle (0.15);
  \draw (1.5,-0.5) circle (3);
  \draw (1.8,-0.7) circle (3.5);
  \filldraw[fill=black!15!white, draw=black] (7,0) circle (1);
  \filldraw[fill=black!15!white, draw=black] (4,-1) circle (0.3);
  \filldraw[fill=black!15!white, draw=black] (4,-3) circle (0.2);
  \filldraw[fill=black!15!white, draw=black](1.5,-2) circle (0.5);
  \draw (1.5,-2) circle (0.7);
  \draw (1.5,-2) circle (0.9);
  \draw (1.5,-2) circle (1.1);
\end{tikzpicture}
\caption{A pattern of circles on $\Sigma$: the innermost disks are colored in grey.}
\label{fig:circles}
\end{center}
\end{figure}
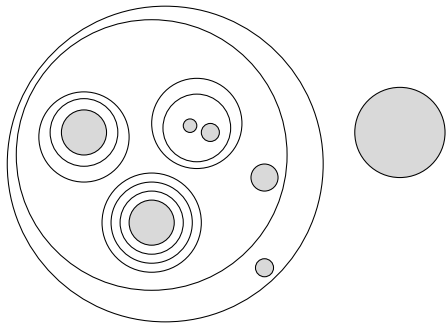

\begin{figure}[h]
\begin{center}
\begin{tikzpicture}[scale=1]
  \draw[very thick] (0,0) ellipse (3 and 2);
  \draw[very thick] (0.7,3.22) arc (140:220:5);
  \draw[very thin] (0,-2) arc (-90:90 :0.2 and 2);
  \draw[very thin, dashed] (0,2) arc (90:270:0.2 and 2);

  \draw (0.3,3) node[above]{$\Sigma$};
  \draw (-3.2,0) node[left]{$s$};
  \draw (-1.9,0) node[right]{$D_1$};
  \draw (0.1,1.2) node[right]{$C$};
  \draw (1.9,0) node[left]{$D_2$};

  \draw[dashed] (8,0) ellipse (3.1 and 2.1);
  \draw[very thick] (9,3.22) arc (140:220:5);
  \draw[very thick] (7.3,1.8) arc (90:270:2.15 and 1.8);
  \draw[very thick] (7.3,0) ellipse (0.15 and 1.8);
  \draw[very thick] (8.7,-1.8) arc (-90:90:2.15 and 1.8);
  \draw[very thick] (8.7,0) ellipse (0.15 and 1.8);
  \draw (8.7,3) node[above]{$\Sigma$};
  \draw (4.8,0) node[left]{$s_1$};
  \draw (11.2,0) node[right]{$s_2$};
  \draw (6,0) node[right]{$D_1$};
  \draw (7.55,0) node[left]{$D$};
  \draw (10,0) node[left]{$D_2$};
  \draw (8.43,0) node[right]{$D$};
\end{tikzpicture}
\caption{The spheres $s$ and $\Sigma$ intersect along a circle $C$. This circle $C$ bounds an innermost disk $D$ on $\Sigma$, and cuts $s$ into two complementary disks $D_1$ and $D_2$. The surgery gives rise to two children $s_1$ and $s_2$, both of them made of a copy of $D$ and a copy of one complementary disk $D_i$.}
\label{fig:surgery}
\end{center}
\end{figure}
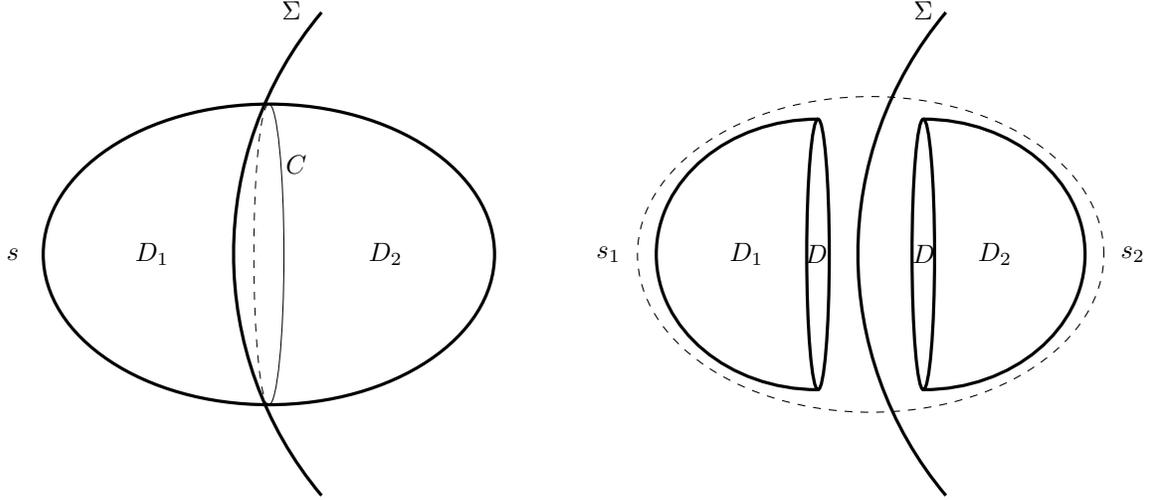

The intersection circles between $S$ and $\Sigma$ define a pattern of circles on $\Sigma$, each of which bounds two disks on $\Sigma$ -- see Figure~\ref{fig:circles}. Choose an innermost disk $D$ in this collection, i.e. the disk $D$ contains no other disk in this pattern, and let $C$ be its boundary circle. The sphere $s\in S$ containing $C$ is thus the union of two disks $D_1$ and $D_2$ which intersect along $C$. Performing surgery on $S$ along $D$ consists of replacing the sphere $s$ by two disjoint spheres $s_1$ and $s_2$ that do not intersect $s$, the sphere $s_1$ being the union of a parallel copy of $D_1$ and a parallel copy of $D$, and $s_2$ being the union of a parallel copy of $D_2$ and a parallel copy of $D$ -- see Figure~\ref{fig:surgery}. We then identify parallel spheres in $S\smallsetminus\{s\}\cup\{s_1,s_2\}$
to get a new sphere system $S'$. We say that $S'$ is obtained from $S$ by performing a single surgery step on $S$ with respect to $\Sigma$. The spheres we obtain from a sphere $s\in S$ after a single surgery step are called the \emph{children} of $s$ (in the case where no surgery is performed on a sphere $s$, then $s$ is its own child).

Given two sphere systems $S$ and $\Sigma$, a \emph{surgery path} from $S$ to $\Sigma$ is a finite sequence $S=S_0,\dots,S_K=\Sigma$ such that for all $i\in \{0,\dots,K-2\}$, the sphere system $S_{i+1}$ is obtained from $S_i$ by performing a single surgery step on $S_i$ with respect to $\Sigma$, and $K$ is the smallest integer such that $S_{K-1}$ does not intersect $\Sigma=S_K$. 
Note in particular that for all $i\in \{0,\dots,K-1\}$, the sphere systems $S_i$ and $S_{i+1}$ do not intersect, so $d(S_i,S_{i+1})\le 2$ (as they are both contained in their union $S_i\cup S_{i+1}$). 

To define a surgery path, after performing a single surgery step, one has to put the new spheres in normal form with respect to $\Sigma$ before continuing to perform other surgeries. The following lemma is classical and can be found for example in \cite{HV96} or in \cite[Lemma 7.5]{HOP12}.

\begin{lemma} \label{surg}
Let $S,\Sigma$ be two sphere systems.
\begin{enumerate}[(i)]
\item If $S'$ is a sphere system obtained by performing a single surgery step on $S$ with respect to $\Sigma$, then $i(S',\Sigma)<i(S,\Sigma)$.

\item There exists a surgery path from $S$ to $\Sigma$ in $\mathcal{S}_n^{(1)}$. 
\end{enumerate}
\end{lemma}

\noindent \textit{Proof} :  Part (i) is a direct consequence of the definition of a single surgery step, and of the fact that putting the sphere system you get after surgery in normal form with respect to $\Sigma$ cannot increase the number of intersection circles. Hence, starting from a sphere system $S$, one can only perform a finite number of surgeries before getting a sphere system that has no intersection with $\Sigma$. This implies part (ii) of the lemma.
\hfill $\square$
\\
\\
\indent We notice that Lemma~\ref{surg} shows that $\mathcal{S}_n^{(1)}$ is connected. This is actually the way this is proved in \cite{Hat95} (where Hatcher even proves that the sphere complex is contractible).
\\
\\
\textit{Remark}: Instead of performing a single surgery step to pass from $S_i$ to $S_{i+1}$, we could have asked to perform simultaneously all surgeries along a collection of  innermost disks. The advantage of this definition is that it enables in particular to define a canonical surgery path from a sphere system to another by performing at each step all innermost surgeries. However, one has to be careful to define these paths properly. Indeed, problems occur when performing simultaneous surgeries on the two opposite sides of the same sphere, or when performing surgery along the last intersection circle of a sphere of $\Sigma$. One way to solve these problems is to use Hatcher and Vogtmann's doubling trick \cite{HV96}: start by adding a parallel copy to every sphere in $S$ before performing surgery on $S$ with respect to $\Sigma$. A slight variation on the argument given in the present paper could be used to show that these paths are unparameterized quasi-geodesics in the sphere complex.
\\
\\
\indent A \emph{generalized surgery path} from $S$ to $\Sigma$ is a finite sequence $S=S_0,\dots,S_K=\Sigma$ such that for all $i\in \{0,\dots,K-2\}$, either $S_{i+1}=S_i$, or else the sphere system $S_{i+1}$ is obtained from $S_i$ by performing a single surgery step on $S_i$ with respect to $\Sigma$, and $S_{K-1}$ does not intersect $S_K=\Sigma$. We will usually denote such a path by $(S_i)_{0\le i\le K}$. In particular, for all $i\in \{0,\dots,K-1\}$, we have $d(S_i,S_{i+1})\le 2$. Note that any surgery path is a generalized surgery path, and conversely, any generalized surgery path obviously defines a surgery path by forgetting the indices $i$ for which $S_{i-1}=S_i$. We say that two generalized surgery paths are \emph{equivalent} if the surgery paths they induce are the same. Note that given two equivalent surgery paths, one can pass from one to the other by inserting and deleting waiting times. 

Let $S'\subseteq S$ be a subsystem of $S$, and $S=S_0,\dots,S_K=\Sigma$ be a generalized surgery path. For any $i\in \{0,\dots,K\}$, we define inductively the \emph{descendants} of $S'$ in $S_i$ in the following way:

$\bullet$ The descendants of $S'$ in $S_0=S$ are the spheres in $S'$.

$\bullet$ The descendants of $S'$ in $S_{i+1}$ are the children of the descendants of $S'$ in $S_i$ if some surgery is performed from $S_i$ to $S_{i+1}$, and are the descendants of $S'$ in $S_i(=S_{i+1})$ otherwise.

\indent Since we identify parallel spheres after every surgery step, two subsystems of $S$ can have common descendants in $S_k$. Understanding these common descendants will be central in our proof of the hyperbolicity of the sphere complex.

The following lemma is a variation on Lemma \ref{surg}. It describes the evolution of intersection numbers along surgery paths. 

\begin{lemma} \label{intersections}
Let $S$,$\Sigma$ be two sphere systems, let $s$ be a sphere in $S$ and let $S=S_0,\dots,S_K=\Sigma$ be a generalized surgery path from $S$ to $\Sigma$. For all $t\in \{0,\dots,K\}$, if $s'$ is a descendant of $s$ in $S_t$, then $i(s',\Sigma)\le i(s,\Sigma)$. Furthermore, if some surgery is performed on $s$ before time $t$, then $i(s',\Sigma)< i(s,\Sigma)$.
\hfill $\square$
\end{lemma}

\indent Let $S,S'$ and $\Sigma$ be three sphere systems. Let $S=S_0,\dots,S_K=\Sigma$ be a generalized surgery path from $S$ to $\Sigma$, and let $S'=S'_0,\dots,S'_{K'}=\Sigma$ be a generalized surgery path from $S'$ to $\Sigma$. We say that the path $(S'_i)_i$ \emph{fellow travels the path $(S_j)_j$ after time $t$} if for all $k\in \{t,\dots,K\}$, the sphere systems $S_k$ and $S'_{k+K'-K}$ share a common sphere. When $K=K'$, we will say that the paths $(S'_i)_i$ and $(S_j)_j$ \emph{fellow travel after time $t$}. Note in particular that any surgery path fellow travels itself after time $0$.

\section{Masur and Minsky's criterion and the projection to a surgery path}

\subsection{Masur and Minsky's criterion for hyperbolicity}
\label{sec:MM}

Let $\mathcal{X}$ be a connected simplicial complex, endowed with the simplicial metric. We define a \emph{path} in $\mathcal{X}$ as a finite sequence $\gamma(0),\dots,\gamma(K)$ of vertices of $\mathcal{X}$ such that $d(\gamma(i),\gamma(i+1))\le 2$. The vertex $\gamma(0)$ is called the \emph{origin} of $\gamma$, and the vertex $\gamma(K)$ is called its \emph{endpoint}. 
A set $\Gamma$ of paths is said to be \emph{transitive} if for all vertices $v,w\in \mathcal{X}$, there exists $\gamma\in \Gamma$ such that $\gamma$ has origin $v$ and endpoint $w$.
A \emph{projection} onto a path $\gamma$ of length $K$ is a map $\pi : \mathcal{X}\to \{0,\dots,K\}$. 

In the following definition, given two integers $a$ and $b$, we will denote by $[a,b]$ the set of integers contained between $a$ to $b$, regardless of whether $a\le b$ or $b\le a$. 

\begin{de}
Let $\mathcal{X}$ be a connected simplicial complex equipped with the simplicial metric. Let $\gamma:\{0,\dots,K\}\to \mathcal{X}$ be a path in $\mathcal{X}$, and let $\pi:\mathcal{X}\to \{0,\dots,K\}$ be a map. Let $A\ge 0$, $B>0$ and $C\ge 0$. We say that $\pi$ is

$\bullet$ \emph{$C$-coarsely retracting} if for all $k\in \{0,\dots,K\}$, the diameter of the set $\gamma([k,\pi(\gamma (k))])$ is less than $C$. 

$\bullet$ \emph{$C$-coarsely Lipschitz} if for all vertices $v,w\in \mathcal{X}$ satisfying $d(v,w)\le 1$, the diameter of the set $\gamma ([\pi(v),\pi(w)])$ is less than $C$.

$\bullet$ \emph{$(A,B,C)$-strongly contracting} if for all vertices $v,w\in \mathcal{X}$ which satisfy $d(v,\gamma ([0,K]))\ge A$ and $d(v,w)\le B\cdot d(v,\gamma ([0,K]))$, the diameter of the set $\gamma ([\pi(v),\pi(w)])$ is less than $C$.
\end{de}

The following theorem due to Masur and Minsky gives a criterion for checking that a connected simplicial metric space is hyperbolic.

\begin{theo} (Masur-Minsky \cite[Theorem 2.3]{MM99})
Let $\mathcal{X}$ be a connected simplicial complex equipped with the simplicial metric. Assume that there exist constants $A\ge 0$, $B>0$, $C\ge 0$, a transitive set of paths $\Gamma$ in $\mathcal{X}$ and for each path $\gamma\in \Gamma$ of length $K$, a map $\pi_\gamma:\mathcal{X}\to \{0,\dots,K\}$, such that all $\pi_\gamma$ are $C$-coarsely retracting, $C$-coarsely Lipschitz and $(A,B,C)$-strongly contracting. Then $\mathcal{X}$ is Gromov hyperbolic, and all the paths $\gamma\in\Gamma$ are unparameterized quasi-geodesics with uniform constants.
\end{theo}

Like Handel and Mosher, we use Masur and Minsky's criterion to prove the hyperbolicity of the sphere complex. However, the collection of paths we consider is the collection of all possible surgery paths, while Handel and Mosher use fold paths. (Lemma \ref{surg} ensures that the family of all surgery paths is transitive, so in particular $\mathcal{S}_n^{(1)}$ is connected.) We thus prove the following theorem (the first part is due to Handel and Mosher \cite{HM12}).

\begin{theo}\label{theo:main}
The sphere complex is Gromov hyperbolic. Furthermore, surgery paths are unparameterized quasi-geodesics with respect to some uniform constants.
\end{theo}

\subsection{The projection to a surgery path} \label{sec-proj}

\begin{figure}
\begin{center}
$\xymatrix{
s'=s'_0\ar[r]&s'_1\ar[r]&\dots\ar[r]&s'_{t'}\ar@{.}[dr]\ar[r]&\dots\ar[r]\ar@{.}[dr]&\Sigma\ar@{.}[dr]&\\
&&&&&&\\
S=S_0\ar[r]&S_1\ar[r]&\dots\ar[r]&S_t\ar[r]\ar@{.}[ur]\ar@{=}[d]&\dots\ar[r]\ar@{.}[ur]&\Sigma\ar@{.}[ur]&\\
&&&\gamma(k)&&&}$
\end{center}
\caption{The projection to a surgery path. The horizontal lines represent the generalized surgery paths from $S$ to $\Sigma$ and from a sphere $s'\in S'$ to $\Sigma$. When two sphere systems are joined by a dotted line, then they are either equal or separated by an edge in $\mathcal{S}_n^{(1)}$ (a ``peak'' between $S_k$ and $s'_{k'}$ indicates a common subsystem shared by $S_k$ and $s'_{k'}$).}
\label{proj}
\end{figure}
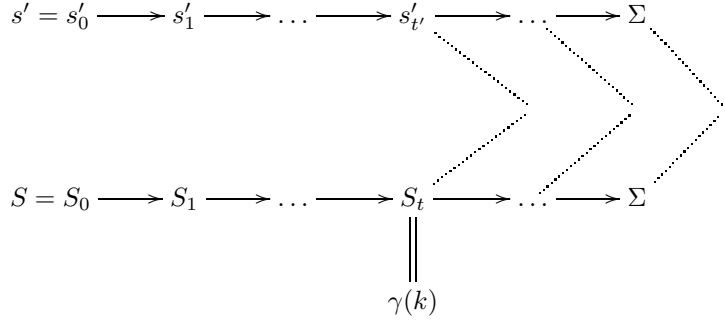

We now have to define the projection to a surgery path in $\mathcal{S}_n^{(1)}$. Let $S,\Sigma$ be two sphere systems, and let $\gamma$ be a surgery path from $S$ to $\Sigma$. Let $S'$ be a sphere system. We define the {\em projection $\pi(S')$ of $S'$ to $\gamma$} to be the smallest integer $k$ such that for each sphere $s'\in S'$, there exist a generalized surgery path $s'=s'_0,\dots,s'_{K'}=\Sigma$ from $s'$ to $\Sigma$, a generalized surgery path $S=S_0,\dots,S_K=\Sigma$ from $S$ to $\Sigma$ equivalent to $\gamma$, and an integer $t$ such that $S_t=\gamma(k)$ and the path $(s'_i)_i$ fellow travels the path $(S_i)_i$ after time $t$ -- see Figure \ref{proj}. 
\\

\noindent \textit{Remark}: Having in mind that surgery paths should resemble folding paths in the reverse direction, our definition of the projection may be seen dual to Handel and Mosher's \cite[Section 4.1]{HM12}. However, the reader might wonder why we are looking at surgery paths starting from all spheres in $S'$ in our definition, and not simply to surgery paths starting from the sphere system $S'$ itself. Actually, using techniques developed in this paper -- combing (Section~\ref{sec-combing}) and understanding common descendants (Section~\ref{sec-son}) -- one can show that these two definitions of the projection are coarsely equal (the proof of one direction is made in Corollary~\ref{cor}, the other direction follows from the arguments in the proof of Proposition \ref{contrW2}). The advantage of our definition is that it will help us for the proof of the coarse retraction property (Proposition~\ref{retr}).
\\
\\
\indent The purpose of the following sections is to prove that the projections to surgery paths are uniformly coarsely retracting (Section \ref{sec-retraction}), coarsely Lipschitz (Section \ref{sec-Lipschitz}) and strongly contracting (Section \ref{sec-contraction}). This finishes the proof of our main result (Theorem~\ref{theo:main}).

\section{Combing} \label{sec-combing}

\begin{de}
Let $k,K\in\mathbb{N}$. A \emph{combing diagram} is a collection $(S_i^j)_{0\le i \le K, 0\le j\le k}$ such that 
\begin{enumerate}[(i)]
\item for all $j\in \{0,\dots,k\}$, the path $(S_i^j)_{0\le i\le K}$ is a generalized surgery path,

\item the sphere systems $S_K^j$ are all equal,

\item for all $j,j'\in\{0,\dots,k\}$ satisfying $|j'-j|=1$ and for all $i\in \{0,\dots,K\}$, if $S_0^j\subseteq S_0^{j'}$, then $S_i^j$ is the sphere system consisting of the descendants of $S_0^j$ in $S_i^{j'}$.
\end{enumerate}
\end{de}

\begin{figure}
\begin{center}
$\xymatrix{
S_0^4\ar[r]\ar@{.}[dr]&S_1^4\ar[r]\ar@{.}[dr]&\ar@{.}[dr]&\dots\ar[r]&\Sigma\ar@{.}[dr]&\\
&S^3_0\ar[r]&S^3_1\ar[r]&&\dots\ar[r]&\Sigma\\
S_0^2\ar[r]\ar@{.}[dr]\ar@{.}[ur]&S_1^2\ar[r]\ar@{.}[dr]\ar@{.}[ur]&\ar@{.}[dr]\ar@{.}[ur]&\dots\ar[r]&\Sigma\ar@{.}[dr]\ar@{.}[ur]&\\
&S^1_0\ar[r]&S_1^1\ar[r]&&\dots\ar[r]&\Sigma\\
S_0^0\ar[r]\ar@{.}[ur]&S_1^0\ar[r]\ar@{.}[ur]&\ar@{.}[ur]&\dots\ar[r]&\Sigma\ar@{.}[ur]&\\}$
\caption{A combing diagram. The horizontal lines represent generalized surgery paths. When we draw a dotted line between two sphere systems, then these sphere systems are either equal or separated by an edge in $\mathcal{S}_n^{(1)}$.}
\label{combing}
\end{center}
\end{figure}
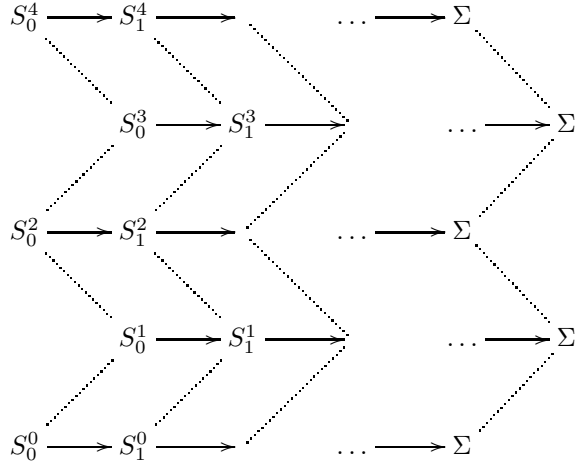

We draw on Figure \ref{combing} such a combing diagram (with $k=4$). 

The goal of this section is, given a zig-zag path between two sphere systems $S$ and $S'$ and a generalized surgery path from $S$ to another sphere system $\Sigma$, to build a combing diagram where the left column is the given zig-zag path and where the bottom row is a generalized surgery path from $S$ to $\Sigma$ equivalent to the given generalized surgery path. This is analogous to \cite[Section 4]{HM12}.

\begin{prop} \label{comb}
Let $S,S'$ and $\Sigma$ be three sphere systems, let $S=S^{0},\dots,S^{k}=S'$ be a zig-zag path in $\mathcal{S}_n^{(1)}$ from $S$ to $S'$, and let $\gamma$ be a generalized surgery path from $S$ to $\Sigma$. Then there exists a combing diagram $(S_i^j)_{0\le i\le K, 0\le j\le k}$ such that

$\bullet$ for all $j\in \{0,\dots,k\}$, we have $S_0^j=S^j$, and

$\bullet$ the generalized surgery path $(S^0_i)_{0\le i\le K}$ is equivalent to $\gamma$.
\end{prop}

\noindent \textit{Proof}: The proof uses two tools, which are analogues for surgery paths of the combings by collapse and expansion of Handel and Mosher's proof \cite[Proposition 4.3 and Proposition 4.4]{HM12}. Note that in this setting, the ``combing by expansion'' step requires less effort than in the original proof.
\\
\\
\textit{``Combing by collapse''}: Let $S$ and $\Sigma$ be two sphere systems, and let $\gamma$ be a generalized surgery path from $S$ to $\Sigma$. Let $S'\subseteq S$ be a nonempty subsystem of $S$. Then $\gamma$ obviously induces a generalized surgery path 
 $\gamma'$
from $S'$ to $\Sigma$ such that for all $i\in \{0,\dots,K\}$, the sphere system $\gamma'(i)$ is a subsystem of the sphere system $\gamma(i)$, by only looking at surgeries made on $S'$ and its descendants, and ``waiting'' at each time when no such surgery is performed.
\\
\\
\textit{``Combing by expansion''}: Suppose now that we are given a generalized surgery path $\gamma'$ from $S'\subseteq S$ to $\Sigma$. We define inductively a generalized surgery path $\gamma$ from $S$ to $\Sigma$ and a generalized surgery path $\widetilde{\gamma}'$ equivalent to $\gamma'$ in the following way. Assume that we have defined $\gamma(i)$ for some $i\in \{0,\dots,K\}$. 

-- If $\gamma'(i+1)=\gamma'(i)$, then we set $\gamma(i+1)=\gamma(i)$, and we do not modify $\widetilde{\gamma}'$.

-- If we use a surgery disk $D$ to pass from $\gamma'(i)$ to $\gamma'(i+1)$, and if $D$ is also innermost when considered in the intersection pattern on $\Sigma$ between $\gamma(i)$ and $\Sigma$ (which contains all the intersection circles between $\gamma'(i)$ and $\Sigma$ because when $\Sigma$ and $\gamma(i)$ are in normal form, then $\Sigma$ and the subsystem $\gamma'(i)$ are also in normal form by the discussion in Section \ref{sec-NF}), then perform the surgery on $\gamma(i)$ to obtain $\gamma(i+1)$. We do not modify $\widetilde{\gamma}'$.

-- If we use a surgery disk $D$ to pass from $\gamma'(i)$ to $\gamma'(i+1)$, and if $D$ is not innermost when considered in the intersection pattern on $\Sigma$ between $\gamma(i)$ and $\Sigma$, then we choose one of the innermost disks contained in $D$. Let $\gamma(i+1)$ be the sphere system obtained from $\gamma(i)$ by performing a single surgery step along this disk. In order to define $\widetilde{\gamma}'$, modify $\gamma'$ by inserting the necessary waiting time between $\gamma'(i)$ and $\gamma'(i+1)$. 

-- If $\gamma'(i+1)=\Sigma$ and $\gamma'(i)\ne\Sigma$, there are again two cases to consider. If $\gamma(i)$ does not intersect $\Sigma$, then we define $\gamma(i+1)=\Sigma$, and we do not modify $\widetilde{\gamma}'$. Otherwise, we choose an innermost disk $D$ in the intersection pattern between $\gamma(i)$ and $\Sigma$, and we let $\gamma(i+1)$ be the sphere system obtained by performing a single surgery step on $\gamma(i)$ along $D$. In order to define $\widetilde{\gamma}'$, modify $\gamma'$ by inserting the necessary waiting time between $\gamma'(i)$ and $\gamma'(i+1)$.

Note that the surgeries performed on $\gamma(i)$ do not affect the intersection circles between $\gamma'(i)$ and $\Sigma$, even after putting the sphere systems in normal form with respect to $\Sigma$. As the intersection number with $\Sigma$ always decreases when we perform some surgery, the process used to define $\gamma$ and $\widetilde{\gamma}'$ terminates. 
\\
\\
The statement of the proposition now follows by iterative applications of ``combing by collapse'' and ``combing by expansion'' (and inserting waiting times when necessary). Note that (iii) in the definition of a combing diagram is automatically guaranteed by our construction.
\hfill $\square$
\\
\\
\noindent \textit{Remark}: When we start with a generalized surgery path of the form $S=S_0,\dots,S_K=\Sigma$ and a zig-zag path from $S=S^0$ to a sphere $S^k$, we will sometimes denote by $(S_i^j)_{0\le i\le K, 0\le j\le k}$ the induced diagram. This is a slight abuse of notation as we should instead consider a generalized surgery path from $S$ to $\Sigma$ equivalent to the path we started with.

\begin{cor} \label{cor}
Let $S$, $S'$ and $\Sigma$ be three sphere systems, and let $\gamma$ be a surgery path from $S$ to $\Sigma$. Let $\pi$ denote the projection to $\gamma$, and let $k=\pi(S')$. Then there exist a generalized surgery path $S'=S'_0,\dots,S'_{K'}=\Sigma$ from $S'$ to $\Sigma$, a generalized surgery path $S=S_0,\dots,S_K=\Sigma$ from $S$ to $\Sigma$ equivalent to $\gamma$, and an integer $t$ such that $S_t=\gamma(k)$ holds and the path $(S'_i)_i$ fellow travels the path $(S_i)_i$ after time $t$.
\end{cor}

\noindent \textit{Proof}: Let $s'$ be a sphere in $S'$. By definition of the projection to a surgery path there exist a generalized surgery path $s'=s'_0,\dots,s'_{K'}=\Sigma$ from $s'$ to $\Sigma$, a generalized surgery path $S=S_0,\dots,S_K=\Sigma$ from $S$ to $\Sigma$ equivalent to $\gamma$, and an integer $t$ such that $S_t=\gamma(k)$ holds and the path $(s'_i)_i$ fellow travels the path $(S_i)_i$ after time $t$. Applying ``combing by expansion'' to the path $(s_i')_i$  yields the desired path from $S'$ to $\Sigma$ (one may furthermore have to add waiting times to the path $(S_i)_i$ to get the desired conclusion).
\hfill $\square$

\section{Coarse retraction} \label{sec-retraction}

We now prove that the projection to a surgery path defined in Section \ref{sec-proj} is coarsely retracting. In our proof, the existence of a sphere $s$ on which no surgery is performed can be seen as an analogue of Handel and Mosher's notion of an ``almost-invariant edge'' \cite[Lemma 5.5]{HM12}. A crucial ingredient of the proof presented here is the fact that the intersection number with the final extremity of a surgery path does not increase along the path (which easily follows from Hatcher's normal form as recalled in Lemma \ref{intersections}). This enables us to avoid some of the technical difficulties encountered in Handel and Mosher's proof.  

\begin{prop} \label{retr}
The projection to any surgery path is $2$-coarsely retracting.
\end{prop}

\noindent \textit{Proof}: Let $S,\Sigma$ be two sphere systems, and let $S=S_0,\dots,S_K=\Sigma$ be a surgery path from $S$ to $\Sigma$. Denote by $\pi$ the projection from $\mathcal{S}_n^{(1)}$ to this surgery path, and let $t\in \{0,\dots,K\}$. The subpath from $S_t$ to $\Sigma$ is a surgery path: thus, by using ``combing by collapse'' for all $s\in S_t$, we see that there exists a generalized surgery path from $s$ to $\Sigma$ which fellow travels $(S_i)_{i\ge t}$. Hence, we have $\pi(S_t)\le t$. 
\\
\\
Let $s\in S_t$ be one of those spheres in $S_t$ that have the fewest intersections with $\Sigma$. By definition of the projection $\pi$, there exist

$\bullet$ a generalized surgery path $s=S'_0,\dots,S'_{K'}=\Sigma$ from $s$ to $\Sigma$, 

$\bullet$ a generalized surgery path $S=\widetilde{S}_0,\dots,\widetilde{S}_{\widetilde{K}}=\Sigma$ equivalent to $(S_i)_i$, and

$\bullet$ an integer $l\in \{0,\dots,\widetilde{K}\}$ 

\noindent such that 

$\bullet$ the sphere system $\widetilde{S}_l$ is equal to $S_{\pi(S_t)}$, and

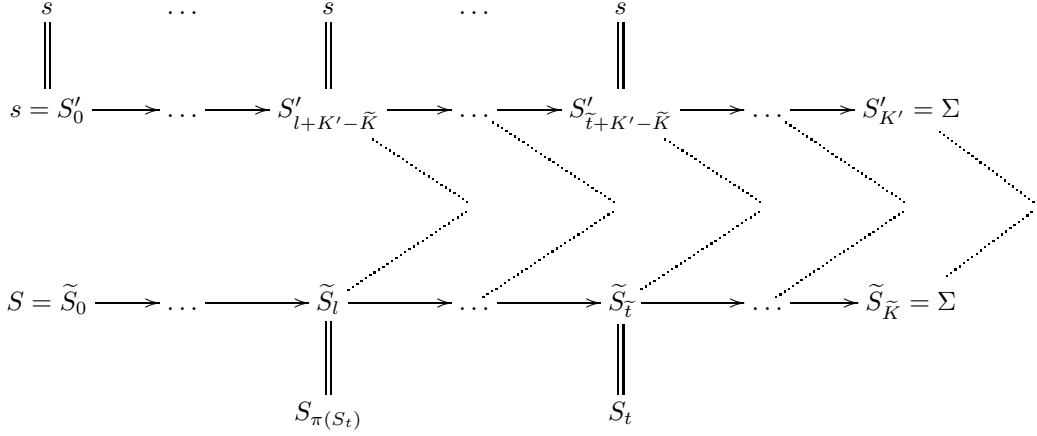
\begin{figure}
\begin{center}
$\xymatrix{
s\ar@{=}[d]&\dots&s\ar@{=}[d]&\dots&s\ar@{=}[d]&&&\\
s=S'_0\ar[r]&\dots\ar[r]&S'_{l+K'-\widetilde{K}}\ar[r]\ar@{.}[dr]&\dots\ar[r]\ar@{.}[dr]&S'_{\widetilde{t}+K'-\widetilde{K}}\ar[r]\ar@{.}[dr]&\dots\ar[r]\ar@{.}[dr]&S'_{K'}=\Sigma\ar@{.}[dr]&\\
&&&&&&&\\
S=\widetilde{S}_0\ar[r]&\dots\ar[r]&\widetilde{S}_l\ar@{.}[ur]\ar[r]\ar@{=}[d]&\dots\ar[r]\ar@{.}[ur]&\widetilde{S}_{\widetilde{t}}\ar[r]\ar@{.}[ur]\ar@{=}[d]&\dots\ar[r]\ar@{.}[ur]&\widetilde{S}_{\widetilde{K}}=\Sigma\ar@{.}[ur]&\\
&&S_{\pi(S_t)}&&S_t&&&}$
\caption{Proof of the coarse retraction property.}
\label{fig-retr}
\end{center}
\end{figure}

$\bullet$ for all $k\in \{l,\dots,\widetilde{K}\}$, the sphere system $S'_{k+K'-\widetilde{K}}$ shares a common sphere with $\widetilde{S}_k$. 
\\
\\
\noindent In particular, let $\widetilde{t}\in \{l,\dots,\widetilde{K}\}$ be such that $S_t=\widetilde{S}_{\widetilde{t}}$, then $S_t=\widetilde{S}_{\widetilde{t}}$ shares a common sphere with $S'_{\widetilde{t}+K'-\widetilde{K}}$ -- see Figure \ref{fig-retr}. This common sphere is a descendant of $s$ in the generalized surgery path $(S'_i)_i$, so that by Lemma~\ref{intersections} it has fewer intersections with $\Sigma$ than $s$ (strictly fewer if some surgery is performed before time $\widetilde{t}+K'-\widetilde{K}$). As $s$ is a sphere in $S_t$ which has minimal number of intersections with $\Sigma$, we derive that the surgery path $(S'_i)_{0\le i\le \widetilde{t}+K'-\widetilde{K}}$ is stationary: for all $i\in\{0,\dots,\widetilde{t}+K'-\widetilde{K}\}$, we have $S'_i =s$. As a result, for all $k\in \{l,\dots,\widetilde{t}\}$, the sphere system $\widetilde{S}_k$ contains $s$, so the diameter of the set $(\widetilde{S}_i)_{l\le i\le \widetilde{t}}$, which is also the diameter of the set $\gamma([\pi(S_t),t])$, is bounded by $2$.

\hfill $\square$

\section{Diameter bounds along surgery paths}

The goal of this section is to provide estimates for the bound of the diameters of surgery paths, in analogy to the bounds exhibited in \cite[Section 5.1]{HM12} for fold paths. Let $S,\Sigma$ be two sphere systems, and let $\gamma$ be a generalized surgery path from $S$ to $\Sigma$. Let $S^1,S^2$ be two nonempty subsystems of $\Sigma$. We say that $S^1$ and $S^2$ \emph{have a common descendant before time t} if the generalized surgery paths induced by $\gamma$ on $S^1$ and $S^2$ have a common (isotopy class of) sphere at time $t$. ``Having no common descendant'' is the analogue of ``having an invariant natural blue-red decomposition'' in Handel and Mosher's proof. As in \cite[Lemma 5.2]{HM12}, we prove in Proposition \ref{son} that as long as two subsystems have no common descendant, the diameter of the path is bounded.  Again, the novelty here is that our approach is based on topological arguments - in particular, Hatcher's normal form for spheres will again play a crucial role in the proof of Proposition \ref{son}. 

\subsection{Complexity of subsystems} \label{sec-complexity}

We start by understanding how the topology of the complement of a sphere system is affected by a single surgery step. In \cite[Lemma 3.1]{HV96}, Hatcher and Vogtmann proved that performing a single surgery on a simple sphere system (i.e. a sphere system $S$ such that the components of $M_n\smallsetminus S$ are simply connected) again yields a simple sphere system (see also \cite[Lemma 3.9]{HH11}). The discussion below may be seen as a generalization of their argument. 

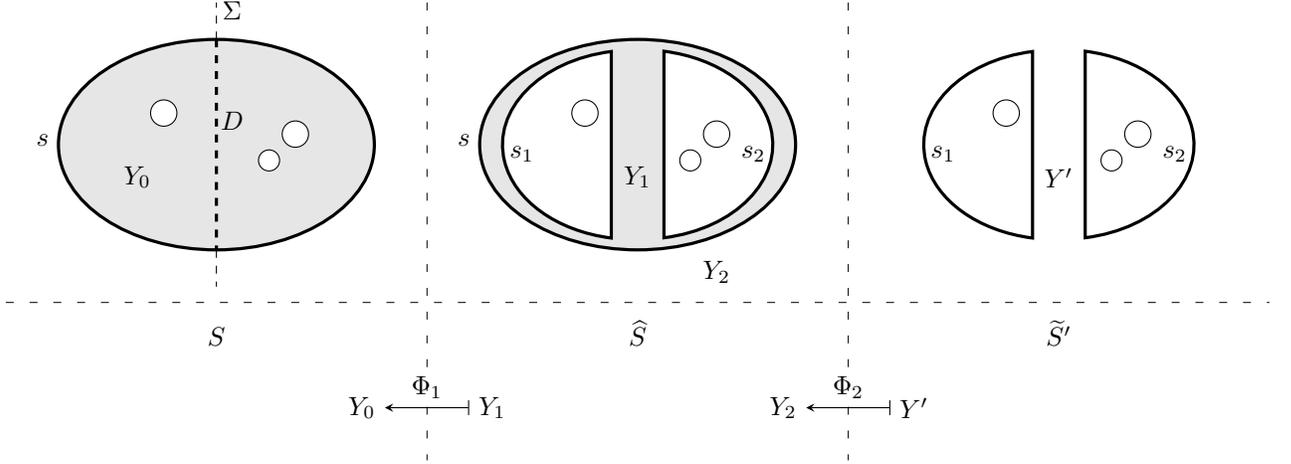
\begin{figure}[h]
\begin{center}
\begin{tikzpicture}[scale=.7]
  \draw[very thick,fill=gray!20] (0,0) ellipse (3 and 2);
  \draw[fill=white] (-1,0.6) circle (.25);
  \draw[fill=white] (1,-0.3) circle (.2);
  \draw[fill=white] (1.5,0.2) circle (.25);
  \draw[very thick, dashed] (0,2) -- (0,-2) ;
  \draw[dashed] (0,2) -- (0,2.7) ;
  \draw[dashed] (0,-2) -- (0,-2.7) ;
  \draw (0.3,2.2) node[above]{$\Sigma$};
  \draw (0.3,0.1) node[above]{$D$};
  \draw (-3.3,-0.2) node[above]{$s$};
  \draw (-1.5,-1) node[above]{$Y_0$};
  
  \draw[very thick,fill=gray!20] (8,0) ellipse (3 and 2);
  \draw[very thick,fill=white] (7.5,1.77) arc (100:260:2.5 and 1.8) -- cycle;
  \draw[very thick,fill=white] (8.5,-1.77) arc (-80:80:2.5 and 1.8) -- cycle;
  \draw[fill=white] (7,0.6) circle (.25);
  \draw[fill=white] (9,-0.3) circle (.2);
  \draw[fill=white] (9.5,0.2) circle (.25);
  \draw (4.7,-0.2) node[above]{$s$};
  \draw (5.8,-0.5) node[above]{$s_1$};
  \draw (10.2,-0.5) node[above]{$s_2$};
  \draw (8,-1) node[above]{$Y_1$};
  \draw (9.5,-2.8) node[above]{$Y_2$};
  
  \draw[very thick,fill=white] (15.5,1.77) arc (100:260:2.5 and 1.8) -- cycle;
  \draw[very thick,fill=white] (16.5,-1.77) arc (-80:80:2.5 and 1.8) -- cycle;
  \draw[fill=white] (15,0.6) circle (.25);
  \draw[fill=white] (17,-0.3) circle (.2);
  \draw[fill=white] (17.5,0.2) circle (.25);
  \draw (13.8,-0.5) node[above]{$s_1$};
  \draw (18.2,-0.5) node[above]{$s_2$};
  \draw (16,-1) node[above]{$Y'$};
  
  \draw[ultra thin,loosely dashed] (-4,-3) -- (20,-3); 
  \draw[ultra thin,loosely dashed] (4,2.7) -- (4,-6); 
  \draw[ultra thin,loosely dashed] (12,2.7) -- (12,-6); 
  \draw (0,-4) node[above]{$S$};
  \draw (8,-4) node[above]{$\widehat{S}$};
  \draw (16,-4) node[above]{$\widetilde{S}'$};
  
  \draw[>=stealth,<-|] (3.2,-5) -- (4.8,-5);
  \draw (4,-5) node[above]{$\Phi_1$};
  \draw (3.2,-5) node[left]{$Y_0$};
  \draw (4.8,-5) node[right]{$Y_1$};
  \draw[>=stealth,<-|] (11.2,-5) -- (12.8,-5);
  \draw (12,-5) node[above]{$\Phi_2$};
  \draw (11.2,-5) node[left]{$Y_2$};
  \draw (12.8,-5) node[right]{$Y'$};

\end{tikzpicture}
\caption{The sets $S$, $\widehat{S}$ and $\widetilde{S}'$.}
\label{fig:components}
\end{center}
\end{figure}

Let $S$, $\Sigma$ be two sphere systems. When performing a single surgery step on $S$ along a disk $D\subseteq\Sigma$, one replaces a sphere $s\in S$ by two spheres $s_1$ and $s_2$. Let $\widehat{S}$ be the set $S\cup\{s_1,s_2\}$ (in which we do not identify parallel spheres), let $\widetilde{S}'$ be the set $S\cup\{s_1,s_2\}\smallsetminus\{s\}$ (again, with no identification of parallel spheres), and let $S'$ be the sphere system derived from $\widetilde{S}'$ by identifying parallel spheres. The sphere system $S'$ is thus a sphere system obtained after a single surgery step on $S$ along $\Sigma$. Denote by $\mathcal{C}$ (resp. $\widehat{\mathcal{C}}$, $\widetilde{\mathcal{C}}'$ and $\mathcal{C}'$) the set of connected components of $M_n\smallsetminus S$ (resp. $M_n\smallsetminus\widehat{S}$, $M_n\smallsetminus\widetilde{S}'$ and $M_n\smallsetminus S'$) -- see Figure \ref{fig:components}. Some of the components in $\widehat{\mathcal{C}}$ and $\widetilde{\mathcal{C}}'$ might be of the form $s'\times [0,1]$ for some sphere $s'\in S'$ ; we call them \emph{pseudo-components}, the other components being \emph{real components}. There is a natural bijection between $\mathcal{C}'$ and the set of real components in $\widetilde{\mathcal{C}}'$, since passing from $\widetilde{S}'$ to $S'$ corresponds to a collapse of pseudo-components. The inclusion $S\subset\widehat{S}$ induces an obvious map $\Phi_1:\widehat{\mathcal{C}}\to\mathcal{C}$. Besides, the components of $M_n\smallsetminus\widetilde{S}'$ are exactly the same as the components of $M_n\smallsetminus\widehat{S}$, except for one component $Y'\in\widetilde{\mathcal{C}}'$. This component $Y'$ is divided into two components $Y_1$,$Y_2\in\widehat{\mathcal{C}}$, where $Y_1$ is the component in $\widehat{S}$ whose boundary consists of the spheres $s$, $s_1$ and $s_2$ (we actually have $Y'=Y_1\cup Y_2\cup s$). We define a map $\Phi_2:\widetilde{\mathcal{C}}'\to\widehat{\mathcal{C}}$ by sending $Y'$ to $Y_2$, and sending every other component to itself. We denote by $\Phi_S:\widetilde{\mathcal{C}}'\to\mathcal{C}$ the composition $\Phi_1\circ\Phi_2$.

Let $Y_0\in\mathcal{C}$ denote the component of $M_n\smallsetminus S$ that contains $D$, and let $Y\in\mathcal{C}\smallsetminus\{Y_0\}$. The preimage $\Phi_S^{-1}(Y)$ has cardinality exactly one. We claim that the component $\Phi_S^{-1}(Y)$ is a real component which has the same rank as $Y$ (the \emph{rank} of a component is defined to be the rank of its fundamental group, which is a free group). This is obvious if $\Phi_S^{-1}(Y)=Y$. The only other possibility is that $\Phi_S^{-1}(Y)=Y'$. This case happens exactly when there are two distinct components in $\mathcal{C}$ that have $s$ as one of their boundary spheres, and when $Y$ is the only component in $\mathcal{C}\smallsetminus \{C_0\}$ having this property. In this case, we have $Y=Y_2$. The claim then follows from the fact that $Y_1$ is simply connected (it is homeomorphic to a ball with punctures), and $Y_2$, whence $Y'$, contains at least three homotopy classes of nontrivial spheres (because $Y_2$ is a real component, as it is equal to the component $Y\in\mathcal{C}$). 

Finally, the component $\Phi_S^{-1}(Y_0)$ either consists of two components whose ranks sum up to the rank of $Y_0$ (by van Kampen's theorem), or it consists of one single component whose rank is one less than the rank of $Y_0$ (this is the case of an HNN extension).

We now define a notion of complexity of a sphere system relative to a proper subsystem. This complexity was introduced by Handel and Mosher in \cite[Section 5.1]{HM12}, and we translate their definition in topological terms. Let $S$ be a sphere system, let $S^1\subsetneq S$ be a proper subsystem, and $S^2:=S\smallsetminus S^1$ be its complementary subsystem. We denote by $C_1(S,S^1)$ the number of connected components of $M_n\smallsetminus S^2$ containing a sphere of $S^1$, and by $C_2(S,S^1)$ the sum of the ranks of these components. The maximal number of spheres in a sphere system is equal to $3n-3$ (this is the dimension of Culler and Vogtmann's unprojectivized Outer space, see \cite{Vog02} for instance). Hence, the maximal number of connected components in the complement of a sphere system is equal to $2n-2$, so we get $0\le C_1(S,S^1)\le 3n-3$. Besides, when gluing two adjacent components along some of their boundary spheres, the rank of the manifold we obtain is greater than the sum of the ranks of the two components we glued. As $M_n$ has rank $n$, iterating this argument shows that $0\le C_2(S,S^1)\le n$. Defining $C(S,S^1):=C_1(S,S^1)+(n-C_2(S,S^1))$, we thus have $0\le C(S,S^1)\le 3n-2$. The following lemma, analogous to \cite[Sublemma 5.3]{HM12}, shows that when performing a single surgery step, the complexity relative to a proper subsystem cannot decrease, provided that the subsystem and its complement subsystem have no common child. The reader should not be surprised that Handel and Mosher's proof naturally translates in the sphere setting, as one single surgery step resembles an elementary unfolding step.

\begin{lemma} \label{complexity}
Let $S$, $\Sigma$ be two sphere systems, and let $S'$ be a sphere system obtained after performing a single surgery step on $S$ along a disk $D\subseteq \Sigma$. Assume that for some proper subsystem $S^1\subsetneq S$, the subsystems $S^1$ and $S^2:=S\smallsetminus S^1$ have no common child in $S'$, and denote by $S'^1$ the subsystem of $S'$ consisting of the children of $S^1$. Then $C(S',S'^1)\ge C(S,S^1)$. Besides, equality occurs if and only if the sphere on which the surgery is performed belongs to $S^1$, or it belongs to $S^2$ and either
\begin{itemize}
\item the component $Y_0$ of $M_n\smallsetminus S^2$ that contains the disk $D$ does not contain any sphere in $S^1$, or

\item $\Phi_{S^2}^{-1}(Y_0)$ consists of two components, one of which is simply connected and does not contain any sphere in $S'^1$.
 \end{itemize}
\end{lemma}
 
\noindent \textit{Proof} : Performing surgery on a sphere in $S^1$ does not affect the components of $M_n\smallsetminus S^2$, and such a component contains a sphere in $S^1$ if and only if it contains a sphere in $S'^1$ after surgery. We now assume we are performing the surgery on a sphere in $S^2$. Let $Y$ be a component of $M_n\smallsetminus S^2$ that does not contain $D$. If $Y$ contains a sphere $s\in S^1$, then $\Phi_{S^2}^{-1}(Y)$ also contains $s$, which is a sphere in $S'^1$ as no surgery is performed on $s$ when passing from $S$ to $S'$. Besides, we know that the components $Y$ and $\Phi_{S^2}^{-1}(Y)$ have the same rank. 

Now assume that the component $Y_0$ of $M_n\smallsetminus S^2$ containing the disk $D$ contains a sphere in $S^1$. Then $\Phi_{S^2}^{-1}(Y_0)$ consists of one or two components, one of which (at least) contains a sphere $s\in S'^1$. If $\Phi_{S^2}^{-1}(Y_0)$ consists of one single component, then it cannot be a pseudo-component, otherwise the sphere $s$ would have a common child with a sphere in $S^2$. So it is a real component, whose rank is one less than the rank of $Y_0$. In this case, we thus have $C(S',S'^1)>C(S,S^1)$. If there are two components $Y$ and $Y'$ in $\Phi_{S^2}^{-1}(Y_0)$, then the same argument as above shows that at least one of them has to be a real component containing a sphere in $S'^1$. We now assume without loss of generality that $Y$ is such a component. If $Y'$ is also a real component which contains a sphere in $S'^1$, then $C_1(S',S'^1)=C_1(S,S^1)+1$ and $C_2(S',S'^1)=C_2(S,S^1)$, whence $C(S',S'^1)>C(S,S^1)$. If $Y'$ does not contain a sphere in $S'^1$, then $C_1(S',S'^1)=C_1(S,S^1)$, and $C_2(S',S'^1)\le C_2(S,S^1)$, with equality if and only if the rank of $Y'$ is equal to zero, i.e. $Y'$ is simply connected.
\hfill $\square$
\\
\\

\subsection{Common descendants and bound estimates} \label{sec-son}

We start by making two more observations about surgery paths.

1) Let $S, \Sigma$ be two sphere systems. Performing a single surgery step on $S$ with respect to $\Sigma$ to get a sphere system $S'$ does not affect property (i) in the definition of normal form with respect to $\Sigma$, and pushing disks that violate (ii) through $\Sigma$ does not affect property (i) either. So the only thing one has to do to put $S'$ in normal form with respect to $\Sigma$ is to push disks that violate (ii) through $\Sigma$. When doing this operation, we call the corresponding disks on $\Sigma$ \emph{vanishing disks}. This operation can also be seen as a surgery that produces a trivial sphere on one side of $\Sigma$.

2) When defining surgery paths, instead of identifying parallel spheres at each step, we can also keep all the copies of parallel spheres and impose that in the sequel of the process, when performing a surgery along a circle of intersection with one of these copies, we perform simultaneously the surgeries along the corresponding circles of intersection with the other parallel copies. Indeed, the notion of ``corresponding circles'', which correspond to parallel circles on $\Sigma$, makes sense by an extension to arbitrary sphere systems of the uniqueness statement for normal form proved in \cite[Proposition 1.2]{Hat95}, see also \cite[Lemma 7.3]{HOP12}. In \cite[Lemma 7.4]{HOP12}, it is proved that if two circles on $\Sigma$ are ``corresponding circles'' that belong to two parallel copies of a sphere $s$, then all the intersection circles on $\Sigma$ that lie between them also belong to spheres that are parallel to $s$. 

Using the two observations above, one may see the set of intersection disks used to perform surgery along the path (called \emph{surgery disks}) and of vanishing disks as a subset of the disks on $\Sigma$ bounded by intersection circles between $S$ and $\Sigma$. For all $t\in \{0,\dots,K-1\}$, the sphere system $S_t$ can be completely recovered from the collection of surgery and vanishing disks used before time $t$: consider the components of the complement in $S$ of the boundary circles of these disks, and cap them off using the disks on $\Sigma$. As a result, if two surgery paths from $S$ to $\Sigma$ are such that the collection of surgery and vanishing disks used before time $t$ are the same, then the sphere system $S_t$ you get at time $t$ is the same in both paths.  

We are now in situation to prove our bound estimates for surgery paths.

\begin{prop} \label{son}
There exists a constant $C_1>0$ such that the following holds.

\noindent Let $S,\Sigma$ be two sphere systems, and let $\gamma$ be a generalized surgery path from $S$ to $\Sigma$. Assume that there exist two nonempty subsystems $S^1,S^2\subseteq S$ having no common descendant before time $t$. Then the diameter of the set $\gamma([0,t])$ is bounded above by $C_1$. 
\end{prop}

\noindent \textit{Proof}: Without loss of generality, we can assume that $S=S^1\amalg S^2$. Indeed, if $S^1$ and $S^2$ are not disjoint, then they already have a common descendant at time $0$. Besides, the sphere system $\gamma(0)$ is at distance at most $1$ of the sphere system $S^1\cup S^2$, and for all $i\in \{0,\dots,t\}$, the sphere system $\gamma(i)$ is at distance at most $1$ of the sphere system which we get at time $i$ in the generalized surgery path starting from $S^1\cup S^2$ and obtained via ``combing by collapse'', which we denote by $\gamma'$. We notice on the one hand that the diameter of $\gamma([0,t])$ is bounded above by two plus the diameter of $\gamma'([0,t])$, and on the other hand that the subsystems $S^1$ and $S^2$ do not have any common descendant before time $t$ in $\gamma'$. 

Using the discussion before the statement of the proposition and the fact that $S^1$ and $S^2$ have no common descendant before time $t$, we can partition the pattern of surgery and vanishing disks on $\Sigma$ used in the path $\gamma$ before time $t$ into two sets: the disks bounded by intersection circles with the spheres in the descendants of $S^1$ on the one hand, and those bounded by intersection circles with the spheres in the descendants of $S^2$ on the other hand.
We subdivide the interval $\{0,\dots,t\}$ into maximal intervals $I_1=\{0=t_0,\dots,t_1\}$, $I_2=\{t_1+1,\dots,t_2\}$,\dots, $I_k=\{t_{k-1}+1,\dots,t_k=t\}$ such that during each of these intervals, we are never using two nested surgery or vanishing disks, where one is bounded by a circle coming from $S^1$, and the other is bounded by a circle coming from $S^2$. (Notice that during an interval, we may perform surgeries both on the descendants of $S^1$ and on those of $S^2$. However we never use one circle on $S^1$ and one circle on $S^2$ which bound nested disks). 

We start by bounding the number $k$ of such intervals. Without loss of generality, we can assume that for at least $\frac{k-1}{2}$ values of $i\in\{2,\dots,k\}$, the first surgery disk we use in the interval $I_i$ is a disk of intersection with a descendant of $S^2$. We will show that for such values of $i$, we have $C(S_{t_i+1},S^1_{t_i+1})>C(S_{t_i},S^1_{t_i})$ (for $t'\in \{0,\dots, t\}$, we denote by $S_{t'}$ the sphere system $\gamma(t')$, and by $S^1_{t'}$ the subsystem of $S_{t'}$ consisting of the descendants of $S^1$ in $S_{t'}$). Together with Lemma \ref{complexity} and the fact that $C(S,S^1)\ge 0$ and $C(S_{t_k},S^1_{t_k})\le 3n-2$, this will imply that $\frac{k-1}{2}\le 3n-2$, hence $k\le 6n-3$. 

Let $i\in\{2,\dots,k\}$ be such that the first surgery disk we use in the interval $I_{i}$ is a disk of intersection with a descendant of $S^2$, which we denote by $D$. During the interval $I_{i-1}$, we would have got the same sphere system at the end if we had started by performing all surgeries coming from $S^2$ before performing all surgeries coming from $S^1$ (provided we do not identify a descendant of $S^1$ and a descendant of $S^2$ even if they are parallel), so we assume that surgeries have been performed in this order. By definition of the intervals, the component $Y_0$ of $M_n\smallsetminus S_{t_i}^2$ that contains $D$ also contains a sphere in $S_{t_i}^1$. By Lemma \ref{complexity}, the only case we have to rule out is the case where $\Phi_{S^2_{t_i}}^{-1}(Y_0)$ consists of two distinct components, one of which, denoted by $Y$, is simply connected and does not contain any descendant of $S^1$. Let $D'$ be a surgery or vanishing disk coming from $S^1$, contained in $D$ and used during the interval $I_{i-1}$. One of the spheres obtained when cutting along $D'$ (which may be trivial if $D'$ is a vanishing disk) is contained in $Y$, otherwise the sphere we are cutting would violate property (i) of the definition of normal form, as $Y$ is simply connected. If this sphere is nontrivial, then its descendants in $S_{t_i+1}$ are again nontrivial spheres contained in $Y$. If it is trivial, then it is obtained from a single surgery on a sphere $s$ which intersects $Y$ only along $D'$, followed by operations to put its descendants in normal form with respect to $\Sigma$. The other sphere obtained from this surgery is a nontrivial sphere contained in $Y$, and its descendants in $S_{t_i+1}$ are again nontrivial spheres contained in $Y$.

It now remains to bound the diameter of the set $(S_i)_{i\in I_l}$, for all $l\in \{1,\dots,k\}$. Let $i,j\in I_l$, with $i\le j$. According to the discussion at the beginning of this section, when going from $S_i$ to $S_j$, we would end with the same sphere system $S_j$ if we had started by performing all surgeries on $S_i$ coming from $S^1$, and then done all surgeries coming from $S^2$. In this way, all the sphere systems we get in the first part contain a common sphere, namely one of the descendants of $S^2$ in $S_i$, on which no surgery is performed. Similarly, all the spheres we get in the second part contain a common sphere. So $d(S_i,S_j)\le 4$. As this is true for all choices of $i\le j$, the diameter of the set $(S_i)_{i\in I_l}$ is bounded by $4$ for all $l\in \{1,\dots,k\}$. 

In particular, we obtain that the diameter of $\gamma([0,t])$ is bounded above by $24n-12$.
\hfill $\square$

\section{Coarse Lipschitz and strong contraction properties}

The goal of this section is to prove that surgery paths are uniformly coarsely Lipschitz and uniformly strongly contracting. Proposition \ref{son} will be our main topological ingredient. Our proof follows roughly the same idea as that in \cite[Section 6]{HM12}: given a zig-zag path between two sphere systems $S$ and $S'$ and a corresponding combing diagram, we will ``contract'' the diagram to ensure that the extremal paths  rapidly fellow travel. The key observation in our approach is that the existence of common descendants (Section \ref{sec-son}) enables us to contract a combing diagram built out of a W-path (Proposition \ref{contrW}), which will let us prove the coarse Lipschitz property in Section \ref{sec-Lipschitz}. An iterated application of this argument will enable us to contract any combing diagram to prove the strong contraction property in Section \ref{sec-contraction}.

\subsection{Contraction of a W-path and proof of the coarse Lipschitz property} \label{sec-Lipschitz}

The following proposition is the key stone in the proof of the coarse Lipschitz and strong contraction properties.

\begin{prop} \label{contrW}
There exists a constant $C_2>0$ such that the following holds.

\noindent Let $S^{0}\supseteq S^{1}\subseteq S^{2}\supseteq S^{3}\subseteq S^{4}$ be a W-path in $\mathcal{S}_n^{(1)}$, and let $\Sigma$ be a sphere system. Let $(S^j_i)_{0\le i\le K, 0\le j\le 4}$ be a combing diagram whose left column is the given W-path. Then there exists $t\in \{0,\dots,K\}$ such that the diameter of $(S^0_i)_{i\le t}$ is less than $C_2$, and such that the paths $(S^0_i)_i$ and $(S^4_i)_i$ fellow travel after time $t$.
\end{prop}

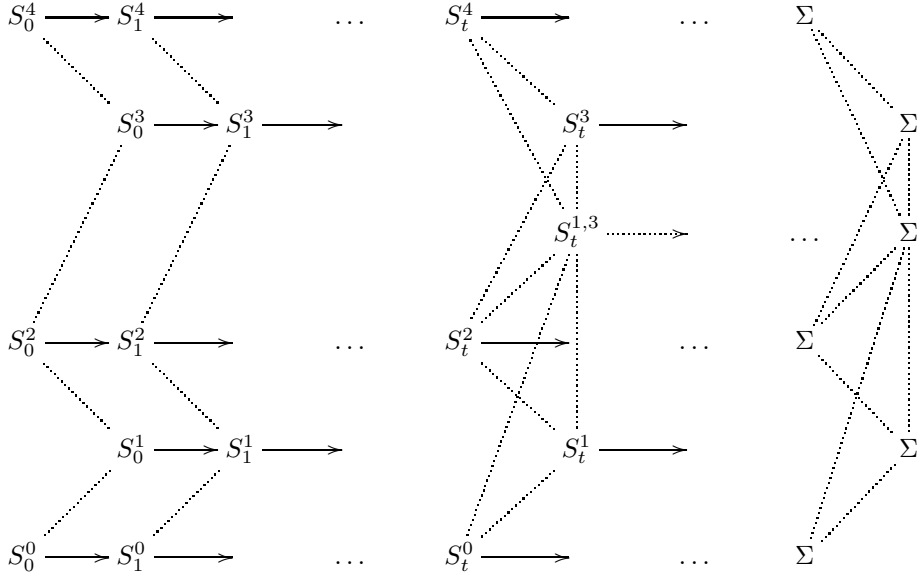
\begin{figure}
\begin{center}
$\xymatrix{
S_0^4\ar[r]\ar@{.}[dr]&S_1^4\ar[r]\ar@{.}[dr]&&\dots&S_t^4\ar[r]\ar@{.}[dr]\ar@{.}[ddr]&&\dots&\Sigma\ar@{.}[dr]\ar@{.}[ddr]&\\
&S_0^3\ar[r]&S_1^3\ar[r]&&&S_t^3\ar[r]\ar@{.}[d]&&&\Sigma\ar@{.}[d]\\
&&&&&S^{1,3}_t\ar@{.>}[r]&&\dots&\Sigma\\
S_0^2\ar[r]\ar@{.}[dr]\ar@{.}[uur]&S_1^2\ar[r]\ar@{.}[dr]\ar@{.}[uur]&&\dots&S_t^2\ar[r]\ar@{.}[dr]\ar@{.}[uur]\ar@{.}[ur]&&\dots&\Sigma\ar@{.}[dr]\ar@{.}[ur]\ar@{.}[uur]&\\
&S_0^1\ar[r]&S_1^1\ar[r]&&&S_t^1\ar[r]\ar@{.}[uu]&&&\Sigma\ar@{.}[uu]\\
S_0^0\ar[r]\ar@{.}[ur]&S_1^0\ar[r]\ar@{.}[ur]&&\dots&S_t^0\ar[r]\ar@{.}[ur]\ar@{.}[uuur]&&\dots&\Sigma\ar@{.}[ur]\ar@{.}[uuur]&}$
\caption{Contracting a W-path.}
\label{fig-W}
\end{center}
\end{figure}

\noindent \textit{Proof}:
Let $C_2:=C_1+4$, where $C_1$ is the constant given by Proposition \ref{son}. Without loss of generality, we may assume that the diameter of the path $(S_i^0)_i$ is greater than $C_2$: otherwise we just set $t=K$. Then, by the triangle inequality, the diameter of the path $(S_i^2)_i$ is greater than $C_1$. Hence, by Proposition~\ref{son}, there exists $t\in \{0,\dots,K\}$ such that $S^{1}_t$ and $S^{3}_t$ have a common descendant in $S^{2}_t$, which we denote by $S^{1,3}_t$. Moreover the diameter of $(S_i^2)_{i\le t}$ is less than $C_1$ (and hence the diameter of $(S^0_i)_{i\le t}$ is less than $C_2$). The surgery paths induced on $S^{1,3}$ by the surgery paths on $S^1$, $S^2$ and $S^3$ in the initial combing diagram are the same, and therefore we get the diagram depicted on Figure \ref{fig-W}. In particular, the paths $(S^0_i)_i$ and $(S^4_j)_j$ fellow travel after time $t$.
\hfill $\square$
\\

We restate the previous proposition in the following way, which is what we will use to prove below the coarse Lipschitz and strong contraction properties.

\begin{prop} \label{contrW2}
There exists a constant $C_3>0$ such that the following holds.\\
Let $S,S'$ and $\Sigma$ be three sphere systems, and let $\gamma$ be a surgery path from $S$ to $\Sigma$. Assume that there exist a generalized surgery path $S=S_0,\dots,S_K=\Sigma$ from $S$ to $\Sigma$ equivalent to $\gamma$, a generalized surgery path $S'=S'_0,\dots,S'_{K'}=\Sigma$ from $S'$ to $\Sigma$, integers $t\in \{0,\dots,K\}$ and $t'\in \{0,\dots,K'\}$, and a W-path from $S_t$ to $S'_{t'}$. Let $k$ be an integer such that $S_t=\gamma(k)$. Let $\pi$ denote the projection of $\mathcal{S}_n^{(1)}$ to $\gamma$. If $\pi(S')\ge k$, then the diameter of $\gamma([k,\pi(S')])$ is bounded above by $C_3$. 
\end{prop}

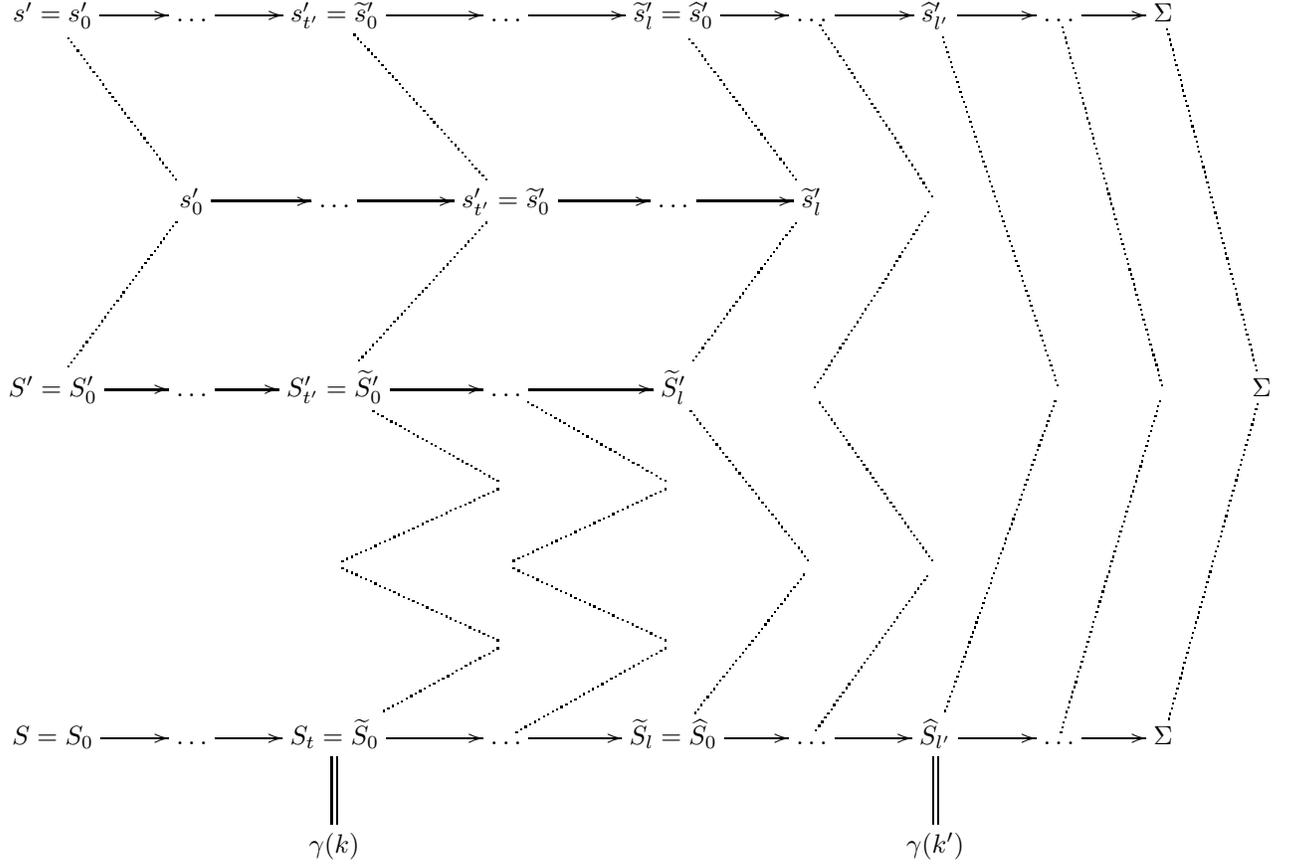
\begin{figure}
\begin{center}
$\xymatrix{
s'=s'_0\ar[r]\ar@{.}[ddr]&\dots\ar[r]&s'_{t'}=\widetilde{s}'_0\ar[r]\ar@{.}[ddr]&\dots\ar[r]&\widetilde{s}'_l=\widehat{s}'_0\ar[r]\ar@{.}[ddr]&\dots\ar[r]\ar@{.}[ddr]&\widehat{s}'_{l'}\ar@{.}[ddddr]\ar[r]&\dots\ar[r]\ar@{.}[ddddr]&\Sigma\ar@{.}[ddddr]\\
&&&&&&&&&\\
&s'_0\ar[r]&\dots\ar[r]&s'_{t'}=\widetilde{s}'_0\ar[r]&\dots\ar[r]&\widetilde{s}'_l&&&&\\
&&&&&&&&&\\
S'=S'_0\ar[r]\ar@{.}[uur]&\dots\ar[r]&S'_{t'}=\widetilde{S}'_0\ar[r]\ar@{.}[dr]\ar@{.}[uur]&\dots\ar[r]\ar@{.}[dr]&\widetilde{S}'_l\ar@{.}[ddr]\ar@{.}[uur]&\ar@{.}[uur]\ar@{.}[ddr]&&&&\Sigma\\
&&&&&&&&&\\
&&\ar@{.}[ur]\ar@{.}[dr]&\ar@{.}[ur]\ar@{.}[dr]&&&&&&\\
&&&&&&&&&\\
S=S_0\ar[r]&\dots\ar[r]&S_t=\widetilde{S}_0\ar[r]\ar@{.}[ur]\ar@{=}[d]&\dots\ar[r]\ar@{.}[ur]&\widetilde{S}_l=\widehat{S}_0\ar[r]\ar@{.}[uur]&\dots\ar[r]\ar@{.}[uur]&\widehat{S}_{l'}\ar@{=}[d]\ar@{.}[uuuur]\ar[r]&\dots\ar[r]\ar@{.}[uuuur]&\Sigma\ar@{.}[uuuur]\\
&&\gamma(k)&&&&\gamma(k')&&&}$
\caption{The diameter of the set $\{\widetilde{S}_i\}_{i\le l}$ is bounded above by $C_2$, as well as the diameter of the set $\{\widehat{S}_i\}_{i\le l'}$.}
\label{big}
\end{center}
\end{figure}

\noindent \textit{Proof}:
Let $C_2$ be the constant provided by Proposition~\ref{contrW}, and let $C_3:=2C_2$. Apply Proposition~\ref{comb} to $(S_i)_i$ and the W-path from $S_t$ to $S'_{t'}$ to obtain a combing diagram, and denote by $(\widetilde{S}_i)_i$ (or $(\widetilde{S}'_i)_i$ resp.) the generalized surgery path from $S_t$ to $\Sigma$ (or from $S'_{t'}$ to $\Sigma$ resp.) obtained by this construction -- see Figure \ref{big}. By Proposition~\ref{contrW}, there exists an integer $l$ such that the paths $(\widetilde{S}'_i)_{i\ge l}$ and $(\widetilde{S}_i)_{i\ge l}$ fellow travel, and the diameter of $(\widetilde{S}_i)_{i\le l}$ is bounded above by $C_2$. The concatenation of $(S'_i)_{i\le t'}$ with $(\widetilde{S}'_i)_i$ is again a generalized surgery path from $S'$ to $\Sigma$, which we denote by $\gamma'$.

Let $s'$ be a sphere in $S'$. Using combing by collapse, we get a generalized surgery path $s'=s'_0,\dots,s'_{t'}=\widetilde{s}'_0,\dots,\widetilde{s}'_l,\dots,\Sigma$ from $s'$ to $\Sigma$ that fellow travels $\gamma'$. In particular, we get a W-path between $\widetilde{S}_l$ and the sphere system $\widetilde{s}'_l$ of the form $\widetilde{S}_l\supseteq\widetilde{S}_l\cap\widetilde{S}'_l\subseteq\widetilde{S}'_l\supseteq\widetilde{s}'_l\subseteq\widetilde{s}'_l$. By the same argument as in the previous paragraph, there exist a generalized surgery path $(\widehat{S}_i)_i$ from $\widetilde{S}_l$ to $\Sigma$ equivalent to a subpath of $\gamma$, a generalized surgery path $(\widehat{s}'_i)_i$ from $\widetilde{s}'_l$ to $\Sigma$, and an integer $l'$ such that the paths $(\widehat{s}'_i)_{i\ge l'}$ and $(\widehat{S}_i)_{i\ge l'}$ fellow travel. Furthermore the diameter of $(\widehat{S}_i)_{i\le l'}$ is bounded above by $C_2$. The concatenation of the path from $s'$ to $\widetilde{s}'_l$ with $(\widehat{s}'_i)_i$ is again a generalized surgery path from $s'$ to $\Sigma$. Let $k'$ be an integer such that $\gamma(k')=\widehat{S}_{l'}$.

Choose $s'\in S'$ to be one of the spheres for which the integer $k'$ one gets by the argument above is the largest. By definition of the projection to a surgery path, we thus have $\pi(S')\le k'$. We also know, by assumption, that $k\le\pi(S')$, and that the diameter of $\gamma([k,k'])$ is bounded above by $2C_2$ (since $\gamma([k,k'])$ is the concatenation of $\{\widetilde{S}_i\}_{i\le l}$ and $\{\widehat{S}_i\}_{i\le l'}$). Therefore, the diameter of $\gamma([k,\pi(S')])$ is bounded above by $C_3$.
\hfill $\square$

\begin{prop} \label{6-4}
There exists a constant $C_3>0$ such that the projection to any surgery path is $C_3$-coarsely Lipschitz.
\end{prop}

\noindent \textit{Proof}: Let $S$ and $\Sigma$ be two sphere systems, and let $\gamma$ be a surgery path from $S$ to $\Sigma$. Let $\pi$ denote the projection from $\mathcal{S}_n^{(1)}$ to $\gamma$. Let $S^1,S^2$ be two sphere systems such that $d(S^1,S^2)=1$. We want to bound above the diameter of $\gamma([\pi(S^1),\pi(S^2)])$. Without loss of generality, we can assume that $\pi(S^1)\le\pi(S^2)$. 

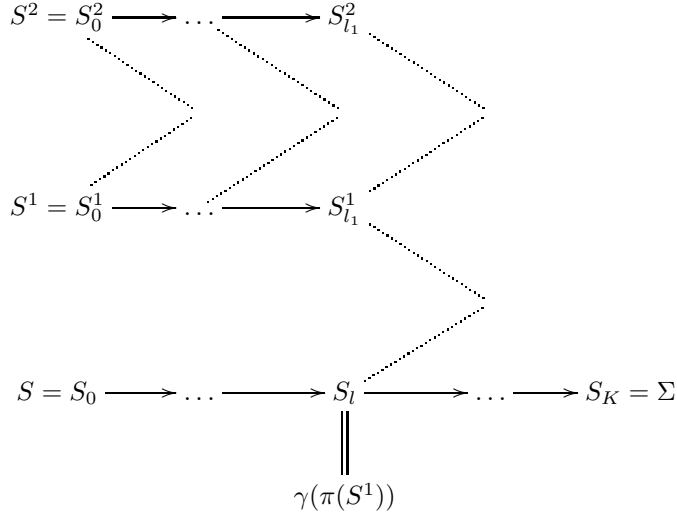
\begin{figure}
\begin{center}
$\xymatrix{
S^2=S^2_0\ar@{.}[dr]\ar[r]&\dots\ar@{.}[dr]\ar[r]&S^2_{l_1}\ar@{.}[dr]&&\\
&&&&\\
S^1=S^1_0\ar@{.}[ur]\ar[r]&\dots\ar@{.}[ur]\ar[r]&S^1_{l_1}\ar@{.}[ur]\ar@{.}[dr]&&\\
&&&&\\
S=S_0\ar[r]&\dots\ar[r]&S_l\ar[r]\ar@{.}[ur]\ar@{=}[d]&\dots\ar[r]&S_K=\Sigma\\
&&\gamma(\pi(S^1))&&}$
\caption{Proof of Proposition \ref{6-4}.}
\label{6.4}
\end{center}
\end{figure}

By Corollary~\ref{cor}, there exist a generalized surgery path $S^1=S^1_0,\dots,S^1_{l_1},\dots,S^1_{K_1}=\Sigma$ from $S$ to $\Sigma$ and a generalized surgery path $S=S_0,\dots,S_l,\dots,S_K=\Sigma$ equivalent to $\gamma$, together with integers $l$ and $l_1$ satisying $K_1-l_1=K-l$ such that $S_l=\gamma(\pi(S^1))$, and $(S^1_i)_i$ fellow travels $(S_j)_j$ after time $l$ -- see Figure \ref{6.4}.

Consider the W-path $S^1\supseteq S^1\cap S^2\subseteq S^2\supseteq S^2\subseteq S^2$ in $\mathcal{S}_n^{(1)}$ (where $S^1\cap S^2$ is either equal to $S^1$ or to $S^2$). Applying Proposition~\ref{comb} to this path and to the generalized surgery path from $S^1$ to $\Sigma$, we get generalized surgery paths from $S^1\cap S^2$ to $\Sigma$ and from $S^2$ to $\Sigma$, together with a combing diagram. In particular, this yields a W-path from $S_l$ to $S^2_{l_1}$. By Proposition~\ref{contrW2}, as we assumed that $\pi(S^1)\le\pi(S^2)$, the diameter of $\gamma([\pi(S^1),\pi(S^2)])$ is at most $C_3$. 
\hfill $\square$

\subsection{Contraction of combing diagrams and proof of the strong contraction property} \label{sec-contraction}

\begin{prop} \label{contrZZ}
There exists a constant $C_4>0$ such that the following holds.

\noindent Let $S,S'$ and $\Sigma$ be three sphere systems, with $S\neq S'$. Let $S=S_0,\dots,S_K=\Sigma$ be a generalized surgery path from $S$ to $\Sigma$, and let $S=S^{0},\dots,S^{2^{k+1}}=S'$ be a zig-zag path from $S$ to $S'$ in $\mathcal{S}_n^{(1)}$, where $k\in\mathbb{N}$ satisfies $2^k-3\le d(S,S')\le 2^{k+1}$. Let $(S^j_i)_i$ be a combing diagram given by Proposition~\ref{comb}. If for some integer $i$, the diameter of the set $\{S^0_0,\dots,S^0_i\}$ is greater than $C_4 d(S,S')$, then the generalized surgery paths from $S=S^0$ to $\Sigma$ and from $S'=S^{2^{k+1}}$ to $\Sigma$ fellow travel after time $i$.
\end{prop}

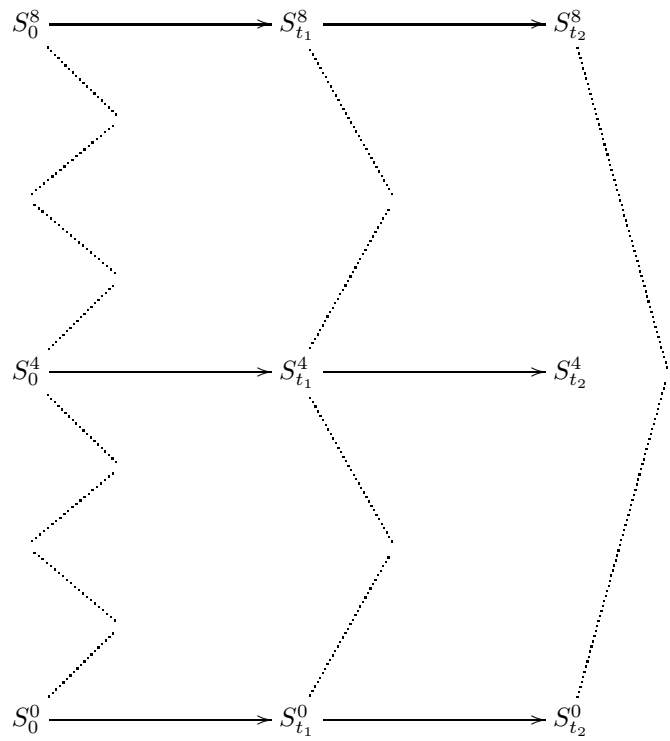
\begin{figure}
\begin{center}
$\xymatrix{
S^8_0\ar[rrr]\ar@{.}[dr]&&&S^8_{t_1}\ar[rrr]\ar@{.}[ddr]&&&S^8_{t_2}\ar@{.}[ddddr]&\\
&&&&&&&\\
\ar@{.}[ur]\ar@{.}[dr]&&&&&&&\\
&&&&&&&\\
S^4_0\ar@{.}[ur]\ar[rrr]\ar@{.}[dr]&&&S^4_{t_1}\ar@{.}[uur]\ar@{.}[ddr]\ar[rrr]&&&S^4_{t_2}&\\
&&&&&&&\\
\ar@{.}[ur]\ar@{.}[dr]&&&&&&&\\
&&&&&&&\\
S^0_0\ar@{.}[ur]\ar[rrr]&&&S^0_{t_1}\ar@{.}[uur]\ar[rrr]&&&S^0_{t_2}\ar@{.}[uuuur]&\\
}$
\caption{Contraction of a zig-zag path made out of two W-paths.}
\label{ZZ2}
\end{center}
\end{figure}

\noindent \textit{Proof}: Let $C_4:=16+8C_2$. The zig-zag path from $S^{0}$ to $S^{2^{k+1}}$ (on the left in Figure \ref{ZZ2}, where we display the situation for $k=2$) is a concatenation of $2^{k-1}$ W-paths. 

If $k=1$, then by Proposition~\ref{contrW} there exists an integer $t_1$ such that the paths $(S^0_i)_i$ and $(S^{2^{k+1}}_j)_j$ fellow travel after time $t_1$, and the diameter of $\{S^0_i\}_{i\le t_1}$ is less than $C_2$.

If $k\ge 2$, then we apply the construction of Proposition~\ref{contrW} to each of the $2^{k-1}$ W-paths to get after some time $t_1$ a diagram made of $2^{k-2}$ bigger W-paths (as well as induced surgery paths starting from these W-paths). We choose $t_1$ to be as small as possible. Then by Proposition~\ref{contrW}, there exists an integer $j_1\in \{0,\dots,2^{k+1}\}$ (which is a multiple of $4$) such that the diameter of $(S^{j_1}_i)_{0\le i\le t_1}$ is smaller than or equal to $C_2$. 

We now iterate the construction above. For $l<k$, we get after some time $t_l$ a diagram made of $2^{k-l-1}$ W-paths (in particular, for all $j,j'\in \{0,\dots,2^{k+1}\}$, we have $d(S_{t_l}^{j},S_{t_l}^{j'})\le 2^{k-l+1}$). We also get an integer $j_l\in \{0,\dots,2^{k+1}\}$ (which is a multiple of $2^{l+1}$) such that the diameter of $(S^{j_l}_i)_{t_{l-1}\le i\le t_l}$ is bounded above by $C_2$.  

We eventually get an integer $t_{k}$ such that the surgery paths from $S$ to $\Sigma$ and from $S'$ to $\Sigma$ fellow travel after time $t_k$, and the diameter of $(S^{0}_i)_{t_{k-1}\le i\le t_k}$ is bounded above by $C_2$. In particular, let $i,i'\in \{0,\dots,t_k\}$ be such that $t_{p-1}\le i\le t_{p}\le\dots\le t_{q}\le i'\le t_{q+1}$ for some integers $p$ and $q$ (we allow the case where $p-1=q$, and we set $t_0=0$). Then we have

\begin{center}
$\begin{array}{rl}
d(S^0_{i},S^0_{i'})&\le d(S^{0}_i,S^{j_{p}}_i)+d(S^{j_{p}}_i,S^{j_{p}}_{t_p})+d(S^{j_p}_{t_p},S^{j_{p+1}}_{t_p})+\dots+d(S^{j_{q}}_{i'},S^{0}_{i'})\\
&\le 2^{k-p+2}+C_2+2^{k-p+1}+C_2+\dots+2^{k-q+1}\\
&\le 2.2^{k+1}+kC_2\\
&\le (2+C_2).2^{k+1}.
\end{array}$
\end{center}

\noindent Using the fact that $2^k-3\le d(S^0,S^{2^{k+1}})$, we thus get

\begin{center}
$\begin{array}{rl}
d(S^0_{i},S^0_{i'})&\le (4+2C_2)d(S^0,S^{2^{k+1}})+12+6C_2\\
&\le (16+8C_2)d(S^0,S^{2^{k+1}}),
\end{array}$
\end{center}

\noindent since we have assumed that $d(S^0,S^{2^{k+1}})\ge 1$. So the diameter of the path $(S^0_i)_{i\le t_k}$ is bounded above by $C_4d(S^0,S^{2^{k+1}})$. 

It follows that if for some integer $i$, the diameter of the set $\{S^0_0,\dots,S^0_i\}$ is greater than $C_4 d(S,S')$, then the generalized surgery paths from $S=S^0$ to $\Sigma$ and from $S'=S^{2^{k+1}}$ to $\Sigma$ fellow travel after time $i$.
\hfill $\square$

\begin{prop}
There exist constants $A\ge 0$, $B>0$ and $C\ge 0$ such that the projection to any surgery path is $(A,B,C)$-strongly contracting.
\end{prop} 

\noindent \textit{Proof}: Let $S,\Sigma$ be two sphere systems, and let $\gamma$ be a surgery path from $S$ to $\Sigma$. Let $S^0$, $S^{2^{k+1}}$ be two distinct sphere systems which satisfy
\begin{center}
$d(S^0,S^{2^{k+1}})\le\frac{1}{C_4+3}d(S^0,\gamma)$.
\end{center}

\noindent By the triangle inequality, we also have

\begin{center}
$\begin{array}{cccc}
d(S^{2^{k+1}},\gamma)&\ge &d(S^0,\gamma)-d(S^0,S^{2^{k+1}})\\
&\ge & (C_4+2) d(S^0,S^{2^{k+1}}).
\end{array}$
\end{center}

\begin{figure}[h]
\begin{center}
\begin{tikzpicture}[scale=0.3]
  \draw[thick] (0,0) node[left] {$S^0$} -- ++(1,-1)
  -- ++(1,1) -- ++(1,-1)
  -- ++(1,1) -- ++(1,-1)
  -- ++(1,1) -- ++(1,-1)
  -- ++(1,1) -- ++(1,-1)
  -- ++(1,1) -- ++(1,-1)
  -- ++(1,1) -- ++(1,-1)
  -- ++(1,1) node[right] {$S^{2^{k+1}}$};
  \draw[thin,dotted] (1,-1) -- ++(1,-1)
  -- ++(1,1) -- ++(1,-1)
  -- ++(1,1) -- ++(1,-1)
  -- ++(1,1) -- ++(1,-1)
  -- ++(1,1) -- ++(1,-1)
  -- ++(1,1) -- ++(1,-1)
  -- ++(1,1);
  \draw[thin,dotted] (2,-2) -- ++(1,-1)
  -- ++(1,1) -- ++(1,-1)
  -- ++(1,1) -- ++(1,-1)
  -- ++(1,1) -- ++(1,-1)
  -- ++(1,1) -- ++(1,-1)
  -- ++(1,1);
  \draw[thin,dotted] (3,-3) -- ++(1,-1)
  -- ++(1,1) -- ++(1,-1)
  -- ++(1,1) -- ++(1,-1)
  -- ++(1,1) -- ++(1,-1)
  -- ++(1,1);
  \draw[thin,dotted] (4,-4) -- ++(1,-1)
  -- ++(1,1) -- ++(1,-1)
  -- ++(1,1) -- ++(1,-1)
  -- ++(1,1);
  \draw[thin,dotted] (5,-5) -- ++(1,-1)
  -- ++(1,1) -- ++(1,-1)
  -- ++(1,1);
  \draw[thick] (0,0) node {$\bullet$} -- ++(12,-12);
  \draw[thick] (14,0) node {$\bullet$} -- ++(-7,-7);
  \draw[thick] (12,-12) -- ++(32,1) node[right] {$\Sigma$} node {$\bullet$};
  \draw[thick] (12,-12) node[below] {$S_l$} -- ++(-16,-0.5) node[left] {$S$} node {$\bullet$};
\end{tikzpicture}
\caption{Schematic representation of the layout.}
\label{fig:strong retraction}
\end{center}
\end{figure}
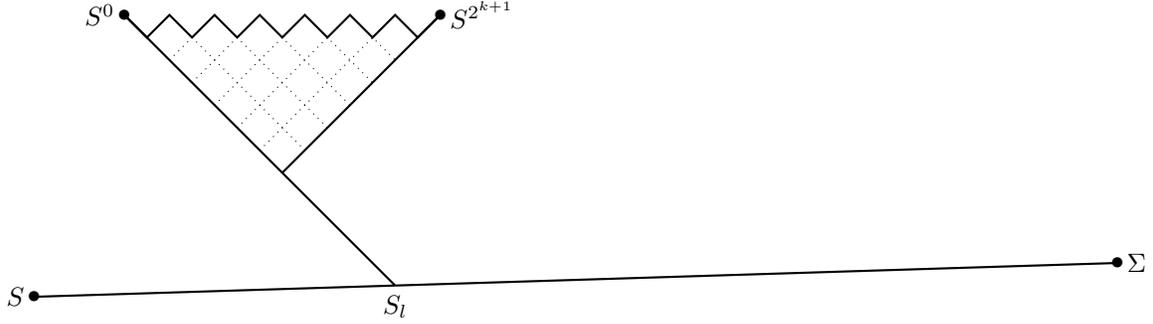

\noindent By Lemma~\ref{ZZ} there exists a zig-zag path of length $2^{k+1}$ from $S^0$ to $S^{2^{k+1}}$, such that $2^k-3\le d(S^0,S^{2^{k+1}})\le 2^{k+1}$. Let $\pi$ denote the projection to $\gamma$, and assume first that $\pi(S^0)\le\pi(S^{2^{k+1}})$. Let $S^0=S^0_0,\dots,S^0_{K_0}=\Sigma$ and $S=S_0,\dots,S_K=\Sigma$ be generalized surgery paths as given by Corollary~\ref{cor}, and let $l$ be an integer such that $S_l=\gamma(\pi(S^0))$, and the path $(S^0_i)_i$ fellow travels the path $(S_i)_i$ after time $l$. Let $(S^j_i)_{0\le i\le K, 0\le j\le 2^{k+1}}$ be a combing diagram given by Proposition~\ref{comb}, obtained from the generalized surgery path $S^0=S^0_0,\dots,S^0_{K_0}=\Sigma$ and the zig-zag path from $S^0$ to $S^{2^{k+1}}$ -- see Figure~\ref{fig:strong retraction}. 

The path $(S^0_i)_{i\le l+K_0-K}$ joins the sphere system $S^0$ to a sphere system at distance at most $2$ from $\gamma$. As $d(S^0,\gamma)\ge (C_4+2) d(S^0,S^{2^{k+1}})$, the diameter of the path $(S^0_i)_{i\le l+K_0-K}$ is thus greater than $C_4 d(S^0,S^{2^{k+1}})$. 
Hence, by Proposition~\ref{contrZZ}, the path $(S^{2^{k+1}}_i)$ fellow travels the path $(S^0_i)$ after time $l+K_0-K$. The sphere systems $S^{2^{k+1}}_{l+K_0-K}$ and $S^0_{l+K_0-K}$ thus share a common sphere, and so do the sphere systems $S^0_{l+K_0-K}$ and $S_l$. So we get a W-path from $S_{l}$ to $S^{2^{k+1}}_{l+K_0-K}$. By Proposition~\ref{contrW2}, the diameter of $\gamma([\pi(S^0),\pi(S^{2^{k+1}})])$ is no more than $C_3$.

In the case when $\pi(S^{2^{k+1}})\le\pi(S^0)$, the proof follows the same scheme: we consider the zig-zag path as a path from $S^{2^{k+1}}$ to $S^0$, and we use that $d(S^0,S^{2^{k+1}})\le\frac{1}{C_4+2}d(S^{2^{k+1}},\gamma)$.

We conclude by setting $A=0$, $B=\frac{1}{C_4+3}$ and $C=C_3$.
\hfill $\square$

\section{Other complexes}\label{sec:other complexes}

In this last section, we mention some consequences of our geometric proof of the hyperbolicity of the sphere complex, concerning the geometry of other related complexes.

\subsection{The free factor complex}\label{sec:free factor}

For $n\in\mathbb{N}$, let $M_{n,1}$ be the manifold obtained from $M_n$ by removing an open ball. Let $\iota:M_{n,1}\to M_n$ denote the inclusion map. We define the \emph{sphere complex} $\mathcal{S}'_{n,1}$ to be the simplicial complex whose vertices are isotopy classes of essential $2$-spheres in $M_{n,1}$ whose $\iota$-image is essential (i.e does not bound a ball) in $M_n$, an $i$-simplex being determined by a system of $i+1$ disjoint, non-isotopic such spheres. We denote by $\mathcal{S}_{n,1}^{(1)}$ the one-skeleton of the first barycentric subdivision of $\mathcal{S}'_{n,1}$ : this is the graph whose vertices are isotopy classes of sphere systems in $M_{n,1}$, two vertices $S$ and $S'$ being joined by an edge if and only if either $S\subsetneq S'$ or $S'\subsetneq S$. Normal form is still well-defined for sphere systems in $\mathcal{S}_{n,1}^{(1)}$ \cite{Hat95}. Surgery paths are also well-defined in $\mathcal{S}_{n,1}^{(1)}$ : performing an elementary surgery step on a sphere $s$ belonging to a sphere system $S\in\mathcal{S}_{n,1}^{(1)}$ yields at least one sphere whose $\iota$-image is essential (it is however possible that the other sphere you get after an elementary surgery has inessential $\iota$-image, in which case it has to be forgotten). 

Let $n\ge 3$. The \emph{free factor complex}, which we denote by $\mathcal{FF}_n$, is defined to be the $\text{Out}(F_n)$-complex whose vertices are conjugacy classes of proper free factors of the free group $F_n$, and $k$-simplices correspond to ascending chains of $k+1$ such free factors (up to global conjugacy), ordered by inclusion.  Hyperbolicity of $\mathcal{FF}_n$ has been proved by Bestvina and Feighn in \cite{BF11}. In \cite{KR12}, Kapovich and Rafi explained how to recover this result from the hyperbolicity of the free splitting complex. We show how our geometric approach gives a new interpretation of Kapovich and Rafi's argument. 

Let $\widehat{\mathcal{S}}_{n}^{(1)}$ be the subcomplex of $\mathcal{S}_n^{(1)}$ consisting of sphere systems whose complement is connected but not simply-connected. There is an obvious inclusion $\widehat{\mathcal{S}}_{n}^{(1)}\subset\mathcal{S}_{n}^{(1)}$. Given a W-path $S^0\supseteq S^1\subseteq S^2\supseteq S^3\subseteq S^4$, we call $S^0$, $S^2$ and $S^4$ the \emph{top vertices}, and $S^1$ and $S^3$ the \emph{bottom vertices}. We extend this terminology to zig-zag paths.

\begin{lemma} \label{lemFF}
Let $n\ge 3$.
\begin{enumerate}[(i)]
\item A sphere system $S\in\mathcal{S}_{n}^{(1)}$ consisting of a single sphere belongs to $\widehat{\mathcal{S}}_{n}^{(1)}$ if and only if the sphere $S$ is non-separating. 

\item Between any two sphere systems $S,S'\in\widehat{\mathcal{S}}_{n}^{(1)}$, there is a zig-zag path in $\mathcal{S}_{n}^{(1)}$ of length at most $4d(S,S')+16$ whose bottom vertices are non-separating spheres (where $d$ is the distance in $\mathcal{S}_{n}^{(1)}$).
\end{enumerate}
\end{lemma}

\noindent \textit{Proof} : Assertion (i) is straightforward. To prove (ii), let $S,S'\in\widehat{\mathcal{S}}_{n}^{(1)}$, and $\gamma$ be a zig-zag path in $\mathcal{S}_{n}^{(1)}$ of the form $S=S^0,\dots,S^{2N}=S'$ whose length is at most $d(S,S')+4$ (the existence of such a path was noticed in the proof of Lemma \ref{ZZ}). We may assume that the bottom vertices of $\gamma$ consist of single spheres, and that the top vertices (except possibly the two extremities) contain at least two spheres. By induction on $k\in\{1,\dots,N-1\}$, we construct disjoint non-separating spheres $s_1^{2k}$, $s_2^{2k}$ and $s'^{2k}$, which either belong to $S^{2k}$ or are disjoint from $S^{2k}$, with the following properties.

\begin{itemize}
\item If $S^{2k-1}$ is non-separating, then $s_1^{2k}=S^{2k-1}$.
\item If $S^{2k+1}$ is non-separating, then $s_2^{2k}=S^{2k+1}$.
\item If $S^{2k-1}$ is separating, then $s_2^{2k-2}$ and $s_1^{2k}$ lie on different sides of $S^{2k-1}$.
\item If $S^{2k}$ contains a non-separating sphere, then $s'^{2k}$ is one of these spheres.
\item If $S^{2k}$ contains no non-separating sphere, then $s_1^{2k}$, $s_2^{2k}$ and $s'^{2k}$ lie in three distinct components of $M_N\smallsetminus S^{2k}$.
\end{itemize}

To do so, we remark that if $\Sigma$ is a system consisting of separating spheres in $M_n$, then each component of $M_n\smallsetminus\Sigma$ contains a non-separating sphere.

By replacing $S^{2k}\supseteq S^{2k+1} \subseteq S^{2k+2}$ by: 
$$S^{2k}\cup s'^{2k}    \supseteq    s'^{2k}    \subseteq    s'^{2k}\cup s^{2k}_2    \supseteq    s^{2k}_2    \subseteq    
s^{2k}_2\cup s^{2k+2}_1    \supseteq    s^{2k+2}_1    \subseteq    s^{2k+2}_1\cup s'^{2k+2}    \supseteq 
s'^{2k+2}    \subseteq    S^{2k+2}\cup s'^{2k+2}$$ 
in $\gamma$ for all $k$, we get a zig-zag path joining $S$ and $S'$ with the desired property.
\qed

\begin{lemma}
For $n\ge 3$, the inclusion $\widehat{\mathcal{S}}_{n}^{(1)}\subset\mathcal{S}_{n}^{(1)}$ is a quasi-isometry.
\end{lemma}

\noindent\textit{Proof} : As every sphere system consisting of one single non-separating sphere belongs to $\widehat{\mathcal{S}}_{n}^{(1)}$, the inclusion map is quasi-surjective (every element of $\mathcal{S}_{n}^{(1)}$ is at distance at most $2$ from an element in $\widehat{\mathcal{S}}_{n}^{(1)}$, because if it does not contain any non-separating sphere, then it is possible to add a non-separating sphere to it). Denote by $d$ the metric in $\mathcal{S}_{n}^{(1)}$, and by $\widehat{d}$ the metric in $\widehat{\mathcal{S}}_{n}^{(1)}$. Let $S,S'\in\widehat{\mathcal{S}}_{n}^{(1)}$. Every geodesic in $\widehat{\mathcal{S}}_{n}^{(1)}$ joining $S$ and $S'$ is a path in $\mathcal{S}_{n}^{(1)}$, so $d(S,S')\le\widehat{d}(S,S')$. Consider a zig-zag path joining $S$ to $S'$ in $\mathcal{S}_{n}^{(1)}$ satisfying the conclusion of Lemma \ref{lemFF}, (ii) : the bottom vertices of this path are single non-separating spheres, and we may also assume that the top vertices (except the two extremities $S$ and $S'$) are systems of two spheres. If a sphere system $S''=\{s_1,s_2\}$ lying in the path consists of two non-separating spheres and does not belong to $\widehat{\mathcal{S}}_n^{(1)}$, then $M_n\smallsetminus S''$ consists of two components $C_1$ and $C_2$, both having nontrivial fundamental group. Let $s$ be a non-separating sphere of $C_1$, then the sphere systems $\{s_1,s\}$ and $\{s_2,s\}$ both belong to $\widehat{\mathcal{S}}_{n}^{(1)}$. We replace $\{s_1\}\subseteq S''\supseteq \{s_2\}$ by a path of length $5$ of the form $\{s_1\}\subseteq\{s_1,s\}\supseteq\{s\}\subseteq\{s_2,s\}\supseteq\{s_2\}$. This yields a path from $S$ to $S'$ of length at most $8d(S,S')+32$ staying in $\widehat{\mathcal{S}}_{n}^{(1)}$. In particular, we have $\widehat{d}\le 8d+32$.
\qed

\begin{theo} (Bestvina-Feighn \cite{BF11})
For all $n\ge 3$, the free factor complex $\mathcal{FF}_n$ is Gromov hyperbolic.
\end{theo}

\noindent\textit{Proof} : Define a map

\begin{displaymath}
\begin{array}{cccc}
\tau:&\widehat{\mathcal{S}}_{n}^{(1)}&\to&\mathcal{FF}_n\\
&S&\mapsto&\pi_1(M_n\smallsetminus S).
\end{array}
\end{displaymath}

\noindent Proposition 3.1 in \cite{HV98} shows that $\tau$ is surjective onto the set of vertices, and $\tau$ is clearly $1$-Lipschitz. 
We remove a small open ball in $M_n$ and we choose a basepoint $p$ on its boundary. For all $S\in\mathcal{S}_n^{(1)}$, we choose a lift $\widetilde{S}$ of $S$ to $\mathcal{S}_{n,1}^{(1)}$. Let $S,S'\in\widehat{\mathcal{S}}_{n}^{(1)}$. Let $\widehat{H}$ denote the conjugacy class of a proper free factor of $F_n$. Assume that $\tau(S)$ and $\tau(S')$ both contain $\widehat{H}$ : there exists $g\in F_n$ such that $\pi_1(M_n\smallsetminus \widetilde{S},p)$ and $g(\pi_1(M_n\smallsetminus \widetilde{S},p))g^{-1}=\pi_1(M_n\smallsetminus g\widetilde{S}'g^{-1},p)$ both contain a free factor $H$ whose conjugacy class is $\widehat{H}$ (where $g\widetilde{S}'g^{-1}$ denotes the image of $\widetilde{S}'$ by a diffeomorphism of $M_{n,1}$ realizing the conjugation by $g$). So by \cite[Lemma 2.2]{HV98}, every element in $H$ is represented by a loop that touches neither $\widetilde{S}$ nor $g\widetilde{S}'g^{-1}$. Hence any sphere system in $\mathcal{S}_{n,1}^{(1)}$ lying on a surgery path joining $\widetilde{S}$ to $g\widetilde{S}'g^{-1}$ has the same property. The inclusion map $M_{n,1}\subset M_n$ yields a projection map from $\mathcal{S}_{n,1}^{(1)}$ to $\mathcal{S}_n^{(1)}$. Notice that the projection to $\mathcal{S}_n^{(1)}$ of any surgery path in $\mathcal{S}_{n,1}^{(1)}$ is a generalized surgery path (performing a surgery in $\mathcal{S}_{n,1}^{(1)}$ for which one of the spheres has to be forgotten because its $\iota$-image is inessential corresponds to a waiting time in the generalized surgery path in $\mathcal{S}_n^{(1)}$). In particular, we get a surgery path $\gamma$ from $S$ to $S'$ in $\mathcal{S}_n^{(1)}$ such that for every sphere system $S''$ in $\gamma$, the image $\tau(S'')$ contains $\widehat{H}$. This implies that the fibers of $\tau$ are quasi-convex. Kapovich and Rafi's criterion \cite[Proposition 2.5]{KR12} then proves the hyperbolicity of $\mathcal{FF}_{n}$.
\qed

\subsection{The arc complex}\label{sec:arc complex}

Our proof of the hyperbolicity of the sphere complex actually shows the following result.

\begin{theo} \label{subcomplex}
Let $\mathcal{X}$ be a subcomplex of $\mathcal{S}'_{n}$ such that
\begin{enumerate}[(i)]
  \item if $\sigma_1$, $\sigma_2$ are two simplices in $\mathcal{X}$ that are faces of a same simplex in $\mathcal{S}'_{n}$, then the simplex spanned by $\sigma_1$ and $\sigma_2$ in $\mathcal{S}'_{n}$ is in $\mathcal{X}$,
  \item if $\sigma_1$, $\sigma_2$ are two simplices in $\mathcal{X}$, and if $\gamma$ is a surgery path from $\sigma_1$ to $\sigma_2$ in $\mathcal{S}'_{n}$, then $\gamma$ is contained in $\mathcal{X}$.
\end{enumerate}
Then $\mathcal{X}$ is Gromov hyperbolic.
\qed
\end{theo}

For $g,s\in\mathbb{N}$, let $S_{g,s}$ be the oriented surface of genus $g$ with $s$ boundary components. An arc on $S_{g,s}$ is said to be \emph{essential} if it is not isotopic, relative to its endpoints, to an arc lying on one boundary component of $S_{g,s}$. The arc complex $\mathcal{A}_{g,s}$, introduced by Harer in \cite{Har85}, is defined to be the simplicial complex whose $k$-simplices correspond to isotopy classes of systems of $k$ essential simple arcs in $S_{g,s}$ with endpoints lying on the boundary of $S_{g,s}$. The manifold $M_{2g+s-1}$ is homeomorphic to the manifold obtained by gluing two copies of $S_{g,s}\times[0,1]$ along their common boundary via the identity map, and the image of an arc under this operation is a sphere. This yields an embedding $i:\mathcal{A}_{g,s}\to \mathcal{S}_{2g+s-1}$. Given two sphere systems $S$ and $\Sigma$ in the image of this map $i$, the pattern of intersection circles between $S$ and $\Sigma$ corresponds to concentric circles on $\Sigma$, and any surgery path from $S$ to $\Sigma$ goes through sphere systems in the image of the map $i$ (each surgery on $S$ corresponds to a surgery on the corresponding arc). In other words, the arc complex $\mathcal{A}_{g,s}$ is a subcomplex of $\mathcal{S}_{2g+s-1}$ satisfying the hypotheses of Theorem \ref{subcomplex}. Theorem \ref{subcomplex} thus yields an alternative proof of the result of Masur and Schleimer that the arc complex is hyperbolic \cite[Theorem 20.3]{MS13}.

\begin{cor} (Masur-Schleimer \cite{MS13})
For all $g,s\in\mathbb{N}$, the arc complex $\mathcal{A}_{g,s}$ is Gromov hyperbolic.
\qed
\end{cor}

\indent It is still open whether the arc complex $\mathcal{A}_{g,s}$ and the sphere complex $\mathcal{S}_{2g+s-1}$ are quasi-isometric in general. Hamenst\"adt and Hensel have proved in \cite[Proposition 4.11]{HH11} that $i$ is a quasi-isometric embedding when $s=1$. The dependence of the quasi-isometry type of the arc complex $\mathcal{A}_{g,s}$ in $g$ and $s$ is also still unknown to us.

\providecommand{\bysame}{\leavevmode\hbox to3em{\hrulefill}\thinspace}

\bigskip
\bigskip

Arnaud Hilion --
LATP -- CMI -- Aix-Marseille Universit\'e -- Technop\^ole Ch\^ateau-Gombert -- 39 rue F. Joliot Curie -- 13453 Marseille Cedex 13 -- France\\
{\tt arnaud.hilion{@}univ-amu.fr}

\bigskip

Camille Horbez --
IRMAR -- Universit\'e de Rennes 1 -- 263 avenue du G\'en\'eral Leclerc, CS 74205 -- 35042 RENNES Cedex -- France\\
{\tt camille.horbez{@}univ-rennes1.fr}

\end{document}